\documentclass[12pt]{amsart}
\usepackage{amsmath,amsthm,amsfonts,amssymb, wrapfig}
\usepackage{graphicx}
\usepackage{wrapfig}
\usepackage{mathrsfs}
\usepackage{enumitem}
\usepackage{wrapfig}
\usepackage{etoolbox}
\usepackage{comment}
\usepackage{xcolor}
\apptocmd{\sloppy}{\hbadness 10000\relax}{}{}
\apptocmd{\sloppy}{\vbadness 10000\relax}{}{}
\usepackage[colorlinks=true, linkcolor=blue, citecolor=blue, urlcolor=blue]{hyperref}
\usepackage[letterpaper,margin=1.1in]{geometry}

\allowdisplaybreaks

\numberwithin{equation}{section}
\theoremstyle{plain}

\newtheorem{theorem}{Theorem}[section]

\newtheorem{corollary}[theorem]{Corollary}

\newtheorem{lemma}[theorem]{Lemma}

\theoremstyle{definition}
\newtheorem{remark}[theorem]{Remark}
\newtheorem{definition}[theorem]{Definition}
\newtheorem{example}[theorem]{Example}

\newtheorem{characterizationproblem}[theorem]{A characterization problem}

\newcommand{\norm}[1]{\left\lVert#1\right\rVert}

\def\SS{\mathbb{S}}

\def\N{\mathbb{N}}

\def\XXint#1#2#3{{\setbox0=\hbox{$#1{#2#3}{\int}$ }
\vcenter{\hbox{$#2#3$ }}\kern-.6\wd0}}

\newcommand{\dist}{\mathop\mathrm{dist}\nolimits}

\newcommand{\loc}{\mathrm{loc}}

\newcommand{\res}{\hbox{ {\vrule height .22cm}{\leaders\hrule\hskip.2cm} }}

\renewcommand{\div}{\operatorname{div}}

\usepackage{alphalph}

\title{Structure of measures for which Ehrhard symmetrization is perimeter non-increasing}
\author{Sean McCurdy and Kuan-Ting Yeh}


\address{Instituto de Matem\'aticas\\ Universidad Nacional Autonoma de M\'exico\\ Cuidad de M\'exico, M\'exico}
\email{sean.mccurdy@im.unam.mx}

\address{Department of Mathematics\\ Purdue University\\
 150 N. University Street, West Lafayette, IN 47907-2067, USA}
 \email{yeh49@purdue.edu}

\begin{document}

\begin{abstract}  
In this paper, we prove that isotropic Gaussian functions are \textit{characterized} by a rearrangement inequality for weighted perimeter in dimensions $n \geq 2$ within the class of non-negative weights in $L^1(\mathbb{R}^n) \cap W^{1,1}_{\loc}(\mathbb{R}^n)$. More specifically, we prove that within this class, generalized Ehrhard symmetrization is perimeter non-increasing for all measurable sets in all directions if and only if the distribution function is an isotropic Gaussian. The class of non-negative $L^1(\mathbb{R}^n) \cap W^{1,1}_{\loc}(\mathbb{R}^n)$-weights is the broadest class in which this problem can be posed for distributional perimeter. One of the main challenges in this paper is handling these weights without imposing any additional structure. Principally, we establish that generalized Ehrhard symmetrization preserves $\mu$-measurability through a novel approximation argument. Additionally, our proof that a rearrangement inequality for weighted perimeter implies that half-spaces are isoperimetric sets is new in the context of generalized Ehrhard symmetrization. Moreover, our version of a variational argument, which had previously appeared in \cite{Rosales14} and \cite{BrockChiacchioMercaldo08}, is carried out under minimal regularity. Finally, we establish some basic but useful results for weighted BV functions with non-negative $L^1(\mathbb{R}^n) \cap W^{1,1}_{\loc}(\mathbb{R}^n)$-weights which may be of independent interest.
\end{abstract}

\maketitle

\vspace{-0.8cm}

\tableofcontents

\section{Introduction}

The main result of this paper establishes that for $n \ge 2$ isotropic Gaussians are the unique distributions that satisfy a symmetrization inequality for their weighted perimeter within the class of non-negative $L^1(\mathbb{R}^n) \cap W^{1, 1}_{\loc}(\mathbb{R}^n)$ weights. Simple examples show that in $n=1$, this uniqueness fails. 

Symmetrization inequalities are a basic tool in the study of isoperimetric problems. In Euclidean space, Steiner symmetrization has played a fundamental role in the study of isoperimetric problems, shape optimization, and rearrangement inequalities  (see, for example, \cite{Steiner1838}\cite{Chlebk2005ThePI}\cite{Krahn25}\cite{PolyaSzego51}\cite{Talenti76} and the hundreds of papers they have inspired).  In Gauss space $(\mathbb{R}^n, \gamma_n)$, the underlying measure is not Lebesgue measure $\mathcal{L}^n$ but $\gamma_n$, defined by 
$$\gamma_n(E)=\frac{1}{(2\pi)^{n/2}}\int_Ee^{-|x|^2/2}\ dx$$
for any measurable set $E\subset \mathbb{R}^n$.
While the Gaussian isoperimetric problem was independently proved by \cite{Borell75_Gaussisoperi} and \cite{CirelsonSudakov74_Gisoperi} without symmetrization, Ehrhard introduced a symmetrization procedure adapted to the underlying measure $\gamma_n$ which generalized Steiner symmetrization.  Ehrhard used this symmetrization procedure to provide an elementary geometric proof of the Gaussian isoperimetric inequality via a Gaussian variation on the Brunn-Minkowski inequality, now called the Ehrhard-Borell inequality \cite{Ehrhard82,Ehrhard83,Ehrhard84}. Since then, Ehrhard's symmetrization procedure (now known as Ehrhard symmetrization) has become an important tool in the study of the Gaussian isoperimetric problems, Gaussian analogs of Faber-Krahn inequalities, and related rearrangement inequalities (see, for example, \cite{Ehrhard84}\cite{CianchiFuscoMaggiPratelli11}\cite{CarlenKerce01}\cite{BrockChiacchioMercaldo12}\cite{DiBlasio03_GaussianTalentiEstimates} among many, many others). Similar statements hold for Schwarz symmetrization, spherical symmetrization, etc., with similarly extensive bibliographies.

The fundamental role of symmetrization inequalities in all these areas leads to a general interest in measures $\mu$, which satisfy a symmetrization inequality for the perimeter.  Roughly speaking, the problem may be phrased as follows:

\begin{characterizationproblem}
    Given a class of measures $\mathscr{W}$ and a generalized symmetrization procedure such that for every $\mu \in \mathscr{W}$ one can define an appropriate notion of weighted $\mu$-perimeter and symmetrization with respect to $\mu$, \textit{characterize} all $\mu \in \mathscr{W}$ such that $\text{Per}_{\mu}(S(E)) \le \text{Per}_{\mu}(E)$
    for all measurable sets $E \subset \mathbb{R}^n$ and all symmetrizations $S$.
\end{characterizationproblem}
For want of a better term, we shall call measures $\mu$ which satisfy their respective perimeter symmetrization inequalities \textit{perimeter symmetrization measures} or \textit{PS measures}.

Solutions to this ``characterization problem" are available in a variety of contexts. In one dimension, Bobkov and Houdr\'e give a characterization of PS measures within the class of probability measures $\mu = f\mathcal{L}^1$ where $f\in C^0((a, b))$ for $f > 0$ on $(a,b)$ (possibly unbounded) for a generalization of Ehrhard symmetrization and the $\varepsilon$-enlargement definition of perimeter (\cite{Bobkov1997} Proposition 13.4).  This builds on Bobkov's previous result characterizing PS measures within the class of log-concave measures $\mu$ on $\mathbb{R}$, under the same choice of symmetrization and perimeter (\cite{Bobkov96} Proposition 2.2). However, in higher dimensions, available results are much more restricted. Though not phrased in this way, Chambers's solution to the Log-convex Density Conjecture proves that, in fact, all smooth radial log-convex measures are PS measures for a generalization of Schwarz symmetrization (appropriately centered) and distributional perimeter \cite{Chambers19}.  Recently, the second author provided a characterization of the equality case in the isoperimetric inequality for \textit{anisotropic} Gaussian functions \cite{Yeh2023}.  As a part of this work, the second author proved that isotropic Gaussian distributions are the unique PS measures within the class of Gaussian distributions on $\mathbb{R}^n$ for a generalized notion of Ehrhard symmetrization and distributional perimeter.

To the best of the authors' knowledge, in all cases in higher dimensions for which the characterization problem has been solved, the class under consideration is very restrictive.  Of course, the geometry of the underlying measure $\mu$ plays a large role, and symmetrization inequalities are expected to be very rigid, at least for $n \ge 2$.  The known examples of PS measures and their symmetrization procedures are expected to be unique.  Explicit computation shows that, for example, the weighted $\mu$-perimeter might not always decrease under generalized Ehrhard symmetrization even if $\mu$ is an anisotropic Gaussian (see \cite{Yeh2023} Example 4.1).  However, expectations and the breakdown of known arguments do not constitute a solution.  

The goal of this paper is to provide a characterization of PS measures for a large class of \textit{finite} measures.

\begin{definition}[The class of measures]
Let $\mathscr{W}(\mathbb{R}^n)$ be the class of measures $\mu$ whose distribution function $\frac{d\mu}{d\mathcal{L}^n} = f$ is non-negative and in $L^1(\mathbb{R}^n)\cap W^{1,1}_{\loc}(\mathbb{R}^n)$. We shall always assume that $f$ is the precise representative (see, for example, \cite{evansgariepy}, Section 1.7.1).
\end{definition}

Because such measures are finite, we will consider a generalized version of Ehrhard symmetrization (see Section \ref{Generalized Ehrhard symmetrization}). As our notion of the perimeter, we shall use the distributional perimeter, that is, the language of weighted BV functions and sets of finite $\mu$-perimeter (see Section \ref{Sets of Finite mu-perimeter and weighted BV functions}).

\begin{theorem}[Main Theorem]\label{Main Theorem}
Let $n \ge 2$ and $\mu \in \mathscr{W}(\mathbb{R}^n)$, as above.  Then, 
\begin{align} \label{e:main theorem symmetrization ineq}
 \operatorname{Per}_{\mu}(S_{\vec{v}}(E))\leq \operatorname{Per}_{\mu}(E)
\end{align}
for all measurable sets $E \subset \mathbb{R}^n$ and all $\vec{v}\in \SS^{n-1}$ if and only if $\mu$ has a distribution function $f(x) = Ce^{-c|x-a|^2}$ for some $0<c<\infty$, $0 \le C <\infty$, and $a \in \mathbb{R}^n$, where $S_{\vec{v}}$ is the generalized $\vec{v}$-Ehrhard symmetrization with respect to $\mu$.
\end{theorem}

The proof that isotropic Gaussian functions satisfy \eqref{e:main theorem symmetrization ineq} is due to Ehrhard \cite{Ehrhard82} if we use the $\varepsilon$-enlargement definition of the perimeter. It is due to \cite{CianchiFuscoMaggiPratelli11} (Theorem 4.3) using the weighted distributional perimeter. As stated above, we shall follow \cite{CianchiFuscoMaggiPratelli11} and use distributional perimeters and the language of weighted BV functions. As such, the contribution of this paper is the following result.

\begin{theorem}[Isotropic Uniqueness]\label{t:main theorem 2}
    Suppose that $n \ge 2$ and $\mu \in \mathscr{W}(\mathbb{R}^n)$.  If generalized $\vec{v}$-Ehrhard symmetrization with respect to $\mu$ is $\mu$-perimeter non-increasing for all measurable sets $E \subset \mathbb{R}^n$ and all $\vec{v} \in \mathbb{S}^{n-1}$, then there exist constants $0<c<\infty$, $0\le C<\infty$, and $a \in \mathbb{R}^n$ such that
    \begin{align*}
        f(x) = Ce^{-c|x-a|^2}.
    \end{align*}
\end{theorem}

\begin{remark}
Theorem \ref{t:main theorem 2}, and therefore Theorem \ref{Main Theorem}, fails in $n=1$. We can construct a probability measure $\mu$ that is a PS measure but does not have the form of a Gaussian distribution (see Example \ref{ex:logconcave}). It is well-known that under reasonable hypotheses, a PS measure on $\mathbb{R}$ must be symmetric (see Lemma \ref{l:Ehrhard implies symmetric}). 
\end{remark}

The proof of Theorem \ref{t:main theorem 2} differs from the proofs of other solutions to the ``characterization problem" in higher dimensions, which rely heavily upon the assumed structure of the measures within their restrictive class.  The second author's solution for general Gaussian functions was accomplished by a detailed study of the relationship between generalized Ehrhard symmetrization and the eigenvectors of the symmetric positive definite matrix $A$ for distributions defined by $e^{-\langle Ax, x\rangle /2}$ \cite{Yeh2023}.  Our solution is most closely related to \cite{Chambers19}, whose proof employs spherical cap symmetrization and a careful study of the ODE one obtains for curves of spherically symmetric perimeter minimizers.  Symmetrization and variational techniques form the heart of our solution as well, but without relying upon smooth distributions or very strong structural assumptions like radial symmetry and log-concavity.

The proof of Theorem \ref{t:main theorem 2} has three steps. Assuming that $\mu$ is a PS measure, we first prove that all half-spaces are perimeter minimizers for the weighted $\mu$-perimeter (Lemma \ref{l:H are minimal}). Our argument follows an intuitive ``see a hole, fill it in" strategy, which is distinct from previous proofs that half-spaces are minimizers for $\mu = \gamma_n$. Previous forms of this argument are available in Euclidean space for Steiner symmetrization (see, for example, \cite{BuragoZalgaller_GeometricInequalitiesBOOK} Lemma 9.4.3). However, in Gauss space, Ehrhard's original geometric proof that half-spaces are $\gamma_n$-perimeter minimizers (\cite{Ehrhard83} Proposition 1.5, though the interested reader should seek out  \cite{lifshits1995gaussian} Chapter 11``Convexity and the Isoperimetric Property" for a detailed English version of the proof) relies essentially upon several key properties of the Gaussian measure $\gamma_n$: it is isotropic, admits a product structure in all orthogonal bases, and is log-concave. Our proof does not rely upon these properties, and is new in the context of generalized Ehrhard symmetrization.  Broadly speaking, Theorem \ref{t:main theorem 2} shows that these properties are a consequence of the symmetrization inequality.  Also, we note that it is the fact that our perimeter minimizers are half-planes which allows us to avoid the more complicated ODE and regularity arguments in \cite{Chambers19}.

The next step is to prove that the distribution function $f$ admits a product structure in all orthogonal bases (Lemma \ref{l:product structure}).  To do this, we use variational techniques like those employed by \cite{Rosales14} in the study of isoperimetric sets for perturbations of log-concave measures and \cite{BrockChiacchioMercaldo12}\cite{BrockChiacchio16} in the study of the structure of measures with foliations by isoperimetric sets. We note, however, that allowing $f$ to vanish introduces several technical challenges.

The last step is to exploit the fact that PS measures in $1$ dimension must be symmetric.  This symmetry, combined with the product structure, forces the distribution function to be radial, which immediately leads to the proof of Theorem \ref{t:main theorem 2}.

It is worth noting that a fundamental challenge in this paper is establishing the measurability of the generalized Ehrhard symmetrization set with respect to $\mu\in \mathscr{W}(\mathbb{R}^n)$. Unlike Steiner and classical Ehrhard symmetrization, its measurability does not immediately follow from the definition and needs to be demonstrated through an approximation argument (see Appendix \hyperref[s:app B measurability]{B}).

The organization of the rest of the paper is as follows.  Section \ref{s:defs and prelims} provides the background for weighted $BV$ functions for weights in $\mathscr{W}(\mathbb{R}^n)$.  This gives us the definition of $\text{Per}_{\mu}$ and some important compactness results, which will be essential to Section \ref{s:half planes are minimizers}.  Section \ref{s:defs and prelims} also gives the definition of generalized Ehrhard symmetrization with respect to $\mu$ and some basic properties of PS measures. Section \ref{s:half planes are minimizers} contains the argument that half-spaces are $\mu$-perimeter minimizers. This argument relies fundamentally upon generalized Ehrhard symmetrization preserving measurability, the proof of which is deferred to Appendix \hyperref[s:app B measurability]{B}. The proof that the distributions of PS measures in $\mathscr{W}(\mathbb{R}^n)$ admit product structures is contained in Section \ref{s:product structure}. Section \ref{s:symm and proof of main theorem} finishes off the proof of Theorem \ref{t:main theorem 2}. 

A weaker version of the main result (assuming the distribution function $f\in C^1(\mathbb{R}^n)\cap W^{1,1}(\mathbb{R}^n)$ is positive) appears in the second author's Ph.D. thesis \cite{Yeh_thesis}. The authors' interest in symmetrization was inspired by recent work in symmetrization within Geometric Measure Theory (see, for example, \cite{Chlebk2005ThePI}, \cite{Barchiesi2013}, \cite{CianchiFuscoMaggiPratelli11}, and \cite{CagnettiColomboDePhilippis17}). It is hoped that this work adds to that interest.

\addtocontents{toc}{\protect\setcounter{tocdepth}{0}}
\section*{Acknowledgments}
The authors would like to thank Tatiana Toro for her feedback on an early version of this paper. A significant part of this work was carried out during S.M.'s postdoctoral fellowship at the National Taiwan Normal University and K.Y.'s studies at the University of Washington, and the authors are grateful to these institutions. Additionally, K.Y. was partially supported by the National Science Foundation Grant DMS--FRG--1853993 during his time at the University of Washington.
\addtocontents{toc}{\protect\setcounter{tocdepth}{1}}

\section{Definitions and preliminaries}\label{s:defs and prelims}

We begin with some basic notation. Let $\chi_{E} $ denote the characteristic function of $E$. We shall use the notation $E \setminus F:= \{x \in E: x\not \in F \}$ and $E \Delta F := (E \setminus F) \cup (E \setminus F)$. For a measure $\mu=f\mathcal{H}^m$ with distribution $f$ and any set $E \subset \mathbb{R}^n$ satisfying $\dim_{\mathcal{H}}(E) = m$, we let $\mu \res E$ be the measure defined by $(\mu \res E)(A) = \int_{E \cap A} f d\mathcal{H}^m$, where $\mathcal{H}^m$ is the $m$-dimensional Hausdorff measure on $\mathbb{R}^n$.  For a vector $\vec{v} \in \mathbb{R}^n$, we denote by $\vec{v}^{\perp}$ the orthogonal complement of $\{\vec{v}\}$, and let $\pi_{\vec{v}^{\perp}}$ be the orthogonal projection onto $\vec{v}^{\perp}$. Note that $\vec{v}^{\perp}=\langle \vec{v}\rangle^\perp=\{\vec{w} \in \mathbb{R}^n: \vec{v} \cdot \vec{w} = 0 \}$. Let $B_r^n(x)$ denote the open ball in $\mathbb{R}^n$ centered at $x$ with radius $r$. When there is no ambiguity, we will simply write $B_r(x)$. The notation $B_r(A):=A+B_r(0)$ denotes the $r$-neighborhood of the set $A$. Finally, $L^p(\mu,U)$ is the classical $L^p$ space with respect to the measure $\mu$ on the open set $U\subset \mathbb{R}^n$.

\subsection{Sets of finite $\mu$-perimeter and weighted BV functions}\label{Sets of Finite mu-perimeter and weighted BV functions}

In this subsection, we recount the basic theory of weighted $BV$ functions for non-negative weights in $W^{1,1}_{\loc}(\mathbb{R}^n)$.  Such weights may be discontinuous, fail to be Doubling and/or Ahlfors regular, and fail to support a weighted Poincar\'{e} inequality.  While there are many excellent books and papers on BV theory, most have concentrated upon weights which satisfy the aforementioned properties.  The authors were unable to find references for weighted BV functions which apply to non-negative $W^{1,1}_{\loc}(\mathbb{R}^n)$ weights in full generality.  Therefore, we include the full proofs of the basic theory in the interest of completeness.

As the reader might expect, the lack of continuity, doubling, Ahlfors regularity, and a weighted Poincar\'{e} inequality puts some limitations on the theory which we are able to develop for such weights. For example the standard proofs of the generalized Gauss-Green formula rely essentially upon the existence of a weighted Poincar\'{e} inequality or a weighted isoperimetric inequality (the two are equivalent for $A_p$ weights by \cite{mazya1973_PoincareIsoperimetric}, but $W^{1,1}$ functions may not be $A_p$ weights).  Fortunately, we do not need these results to prove Theorem \ref{Main Theorem} or Theorem \ref{t:main theorem 2}. The main results we need are weak compactness (Lemma \ref{l:weak compactness}) and the trace operator or Gauss-Green formula for Sobolev functions (Lemma \ref{l:structure theorem}).

\begin{definition}[Weighted $BV$ functions]\label{Weighted_BV}
Let $\mu \in \mathscr{W}(\mathbb{R}^n)$ with a distribution function $f$. We shall say a function $g: \mathbb{R}^n \rightarrow \mathbb{R}$ has {\bf bounded $\mu$-variation} in an open set $U \subset \mathbb{R}^n$ if and only if
\begin{align*}
    \sup \left\{\int_U g\Big(\text{div}(\phi)f +\phi \cdot \nabla f\Big)dx: \phi\in C^1_c(U; \mathbb{R}^n), \|\phi\|_{\infty} \le 1\right\}<\infty.
\end{align*}
Here $\|\cdot\|_\infty$ is the $L^\infty$-norm with respect to the Lebesgue measure on $\mathbb{R}^n$. We shall denote the collection of all functions with bounded $\mu$-variation in $U$ by $BV_{\mu}(U)$.  For unbounded sets $U$, we shall say $g\in BV_{\mu, loc}(U)$ if $g \in BV_{\mu}(V)$ for all open $V \subset \subset U$. If $g \in BV_{\mu}(U)$, then we will use the notation
\begin{align*}
\norm{g}_{BV_\mu(U)} := \sup \left\{\int_{U}g\Big( \text{div}(\phi)f + \phi \cdot \nabla f\Big) dx: \phi \in C^1_c(U; \mathbb{R}^n), \|\phi\|_{\infty} \le 1\right\}.
\end{align*}

For a measurable set $E \subset \mathbb{R}^n$, we shall say that $E$ is a {\bf set of finite $\mu$-perimeter} in an open set $U \subset \mathbb{R}^n$ if and only if $\chi_{E} \in BV_{\mu}(U)$. If $\chi_{E} \in BV_{\mu}(U)$, we shall use the notation
    \begin{align}\label{e:perimeter BV def}
       \operatorname{Per}_{\mu}(E; U):= \norm{\chi_{E}}_{BV_{\mu}(U)}.
    \end{align}
\end{definition}

\begin{lemma}[Riesz Representation Theorem]\label{RRT}
    Let $\mu \in \mathscr{W}(\mathbb{R}^n)$ with distribution function $f$.  For all open sets $U \subset \mathbb{R}^n$ and all $g \in BV_{\mu}(U)$, there exists a Radon measure $\norm{Dg}_\mu$ on $U$ and a $\norm{Dg}_\mu$-measurable function $\sigma: U \rightarrow \mathbb{R}^n$ such that 
    \begin{enumerate}
        \item $|\sigma(x)| = 1$ for $\norm{Dg}_\mu$-a.e. $x \in U$.
        \item For all $\phi \in C^1_c(U; \mathbb{R}^n)$ 
        \begin{align*}
            \int_{U}g\Big(f \operatorname{div}(\phi) + \phi \cdot \nabla f\Big) dx = - \int_U \phi \cdot \sigma d\norm{Dg}_{\mu}.
        \end{align*}
    \end{enumerate}
    Note that $\norm{Dg}_\mu(U) = \norm{g}_{BV_{\mu}(U)}$.
\end{lemma}

The proof follows the standard proof for $f \equiv 1$ in \cite{evansgariepy}, Section 5.1, Theorem 1. We will omit the notation $\mu$ in $\norm{Dg}_\mu$ and simply write $\norm{Dg}$ if there is no confusion.

\begin{lemma}[Lower Semicontinuity]\label{l:lsc}
Let $\mu \in \mathscr{W}(\mathbb{R}^n)$ with distribution function $f$, and let $U \subset \mathbb{R}^n$ be open.  Let $g_k \in BV_{\mu}(U)$ be a sequence that satisfies the following two conditions.
    \begin{enumerate}
        \item $\sup_k \norm{g_k}_{\infty}< \infty$.
        \item There exists a function $g \in L^1_\loc(\mu, U)$ such that \begin{align*}
        g_k \rightarrow g \quad \text{in  } L^1_{\loc}(\mu, U).
    \end{align*}
    \end{enumerate}
    Then $\norm{Dg}(U) \le \liminf_{k \rightarrow \infty}\norm{Dg_k}(U)$.
\end{lemma}

\begin{proof}
First, we note that if $\sup_k \norm{g_k}_{\infty}< \infty$ then $\norm{g\chi_{\{f>0\}}}_{\infty} \le \sup_k \norm{g_k}_{\infty}$. Now, let $0<R<\infty$ be fixed.  Note that since $f \ge 0$ we may write $f = \chi_{\text{spt}(f)}e^h$ for some function $h: \mathbb{R}^n \rightarrow \mathbb{R}$. Therefore, $\nabla f = (\nabla h) f$ and $\nabla h\in L^1_{\loc}(\mu, \mathbb{R}^n)$.  Additionally, observe that for any $\phi \in C^1_c(U; \mathbb{R}^n)$ such that $\norm{\phi}_{\infty} \le 1$, $\text{div}(\phi)$ is bounded. Therefore,
$$\left|\int_{B_R(0) \cap U} (g-g_k)\text{div}(\phi)fdx\right| \le \norm{\text{div}(\phi)}_{\infty}\norm{g-g_k}_{L^1(\mu, B_R(0) \cap U)} \rightarrow 0$$
and for any $K>0$,
\begin{align*}
    \left|\int_{B_R(0)\cap U} (g-g_k)(\phi \cdot \nabla f)dx\right| & = \left|\int_{B_R(0)\cap U} (g-g_k)(\phi \cdot \nabla h)fdx\right|\\
    &
\le \norm{(g-g_k)\chi_{\{f>0\}}}_{\infty}\int_{B_R(0) \cap U \cap \{|\nabla h|>K\}}|\nabla{h}|fdx\\
&\quad+ K \int_{B_R(0)\cap U}|g-g_k|fdx.
\end{align*}
Taking $K$ sufficiently large such that $\norm{\nabla h}_{L^1(\mu, B_R(0) \cap \{|\nabla h|>K\})} \le \varepsilon$, we obtain
\begin{align*}
    \lim_{k \rightarrow \infty}\left| \int_{B_R(0) \cap U} (g-g_k)(\phi \cdot \nabla f)dx\right| \le 2\sup_k \norm{g_k}_{\infty} \varepsilon,
\end{align*}
where we have used $\norm{(g-g_k)\chi_{\{f>0\}}}_{\infty}\leq 2\sup_k \norm{g_k}_{\infty}$. Letting $\varepsilon \rightarrow 0$, we see that 
\begin{align*}
    \int_{B_R(0) \cap U}g(f\text{div}(\phi) + \phi \cdot \nabla f)dx & = \lim_{k \rightarrow \infty} \int_{B_R(0) \cap U}g_k(f\text{div}(\phi) + \phi \cdot \nabla f )dx\\
    & = -\lim_{k \rightarrow \infty} \int_{B_R(0) \cap U}\phi \cdot \sigma_k d\norm{Dg_k} \\
    & \le \liminf_{k \rightarrow \infty} \norm{Dg_k}(B_{R}(0) \cap U),
\end{align*}
where we have applied Lemma \ref{RRT} on $g_k$. Taking the supremum in $\phi$ and letting $R \rightarrow \infty$, we obtain the claim.
\end{proof}

    The next result (Lemma \ref{l:weak approx by smooth}) follows the outline of Evans and Gariepy \cite{evansgariepy}, Section 5.2.2, Theorem 2.  We note that Lemma \ref{l:weak approx by smooth} is essential for proving Weak Compactness of $BV_{\mu}$ in Lemma \ref{l:weak compactness}.  

\begin{lemma}[Weak approximation by smooth functions has bounded weighted variation]\label{l:weak approx by smooth} Let $\mu \in \mathscr{W}(\mathbb{R}^n)$ with distribution function $f$.  If $U \subset \mathbb{R}^n$ is open and $g \in BV_{\mu}(U) \cap L^\infty(U)$, then there is a sequence of smooth functions $\{g_\varepsilon\}_{\varepsilon > 0} \subset BV_{\mu}(U)$ such that 
    \begin{enumerate}
        \item $g_\varepsilon \rightarrow g$ in $L^1_{\loc}(\mu, U)$ as $\varepsilon \rightarrow 0$, and $\norm{g_\varepsilon}_{\infty} \le \norm{g}_{\infty} (1+2\varepsilon)$ for all $\varepsilon>0$.
        \item $\limsup\limits_{\varepsilon \rightarrow 0}\norm{Dg_\varepsilon}(U) \le \norm{Dg}(U) + C(n)\norm{g}_{\infty}\norm{\nabla f}_{L^1(U)}$.
    \end{enumerate}
In fact, if $U$ is bounded, then $g_\varepsilon \rightarrow g$ in $L^1(\mu, U)$.  In the case that $U$ is unbounded, the estimate in (2) may be entirely vacuous if $f \not \in W^{1,1}(\mathbb{R}^n)$.
\end{lemma}

\begin{proof}
    Fix $\varepsilon>0$, and let $m \in \mathbb{N}$.  For $k = 1, 2, ...$ we define 
    \begin{align*}
        U_k:= \left\{x \in U: \dist(x, \partial U) \ge \frac{1}{m+k} \right\} \cap B_{m+k}(0).
    \end{align*}
Now, choose $m$ sufficiently large such that $\norm{Dg}(U \setminus U_1) \le \varepsilon$.  Letting $U_0 = \emptyset$, we inductively define 
\begin{align*}
    V_k = U_{k+1} \setminus \overline{U}_{k-1}
\end{align*}
for $k \in \mathbb{N}$. Let $\xi_k$ be a partition of unity subordinate to $\{V_k\}_k$. Note that we may choose $|D\xi_k| \le C\cdot(m+k)$ for some $C=C(n)>0$ and $\sum_{k}\chi_{V_k}(x) \le 2$.  

Now, let $\eta$ be a standard smooth, radial mollifier with $\text{supp}(\eta) \subset B_1(0)$ and $0 \le \eta\le 1$.  Let $\eta_\varepsilon(x) = \varepsilon^{-n}\eta(x/\varepsilon )$ and choose $\varepsilon_k$ such that the following conditions hold.
\begin{enumerate}
    \item (Containment) $\text{supp}(\eta_{\varepsilon_k} \star \xi_k) \subset V_k$.
    \item (Convolution Estimate I) For all $k \in \mathbb{N}$ 
    \begin{align}\label{e: POU 2}
    \int_U |\eta_{\varepsilon_k} \star (g\xi_k) - g\xi_k|dx \le \frac{\varepsilon}{2^k}.
\end{align}
\item  (Convolution Estimate II) For all $k \in \mathbb{N}$ 
    \begin{align}\label{e:POU 3}
    \int_U |\eta_{\varepsilon_k} \star (fg\nabla\xi_k) - fg\nabla\xi_k|dx \le \frac{\varepsilon}{2^k}.
\end{align} 
\item (Modulus control) For all $k 
\in \mathbb{N}$, if $x \in \{\xi_k \ge \varepsilon/2 \}$ then
\begin{align*}
|\eta_{\varepsilon_k} \star \xi_k(x)| \le |\xi_k(x)|(1+\varepsilon)
\end{align*}
and if $x \in \{\xi_k < \varepsilon/2 \}$ then
\begin{align*}
|\eta_{\varepsilon_k} \star \xi_k(x)| \le \varepsilon.\end{align*}

\end{enumerate}
Now, we define
\begin{align}
    g_\varepsilon := \sum_{k=1}^{\infty}\eta_{\varepsilon_k}\star (g \xi_k).
\end{align}
Note that $g_{\varepsilon} \in C^\infty(U)$ and therefore $g_{\varepsilon} \in BV_{\mu, loc}(U)$ with $\norm{Dg_{\varepsilon}}(V) \le \norm{Dg_{\varepsilon}}_{L^1(\mu, V)}$ for all $V \subset \subset U$. 

To see that $\norm{g_{\varepsilon}}_{\infty} \le \norm{g}_{\infty}(1+2\varepsilon)$, we consider the properties of $\{\xi_k\}_k$ as a partition of unity and recall that $\sum_{k}\chi_{V_k}(x) \le 2$. For any $x \in V_k$, there are three possibilities:
\begin{enumerate}
    \item[(a)] If $x \in \{\xi_k < \varepsilon/2 \}$ then $x \in \{\xi_{k+1} \ge \varepsilon/2 \}$, and $\xi_j(x)=0$ for all $j \not = k, k+1$.
    \item[(b)] If $x \in \{\xi_k > 1- \varepsilon/2 \}$ then $x \in \{\xi_{k+1} < \varepsilon/2 \}$, and $\xi_j(x)=0$ for all $j \not = k, k+1$.
    \item [(c)] If $x \in \{\xi_k \le 1- \varepsilon/2 \}$ then $x \in \{\xi_{k+1} \ge \varepsilon/2 \}$, and $\xi_j(x)=0$ for all $j \not = k, k+1$.
\end{enumerate}
In all three cases, we obtain from condition (4) that $\sum_k |\eta_k \star \xi_k(x)| \le 1+2\varepsilon$. Therefore, because $|\eta_{\varepsilon_k} \star g \xi_k| \le \norm{g}_{\infty} |\eta_{\varepsilon_k} \star \xi_k|$ and $k$ was arbitrary, we see that $$\norm{g_{\varepsilon}}_{\infty} \le \norm{g}_{\infty}(1+2\varepsilon).$$

In order to see that $g_\varepsilon \rightarrow g$ in $L^1_\loc(\mu, U)$, for each $\delta>0$ we let $K=K(\delta)$ be such that $\int_{\{f\ge K\}}fdx \le \delta$.  Then, we may calculate $\norm{g-g_{\varepsilon}}_{L^1(\mu, U)}$ using \eqref{e: POU 2} and the fact that $\sum_k \chi_{V_k} \le 2$.
\begin{align*}
    \norm{g-g_{\varepsilon}}_{L^1(\mu, U)} & \le
    \int_{U}\sum_{k}|\eta_{\varepsilon_k}\star(g\xi_k) - g\xi_k|fdx\\
    & \le \int_{U \cap \{f < K \}}\sum_{k}|\eta_{\varepsilon_k}\star(g\xi_k) - g\xi_k|fdx\\
    & \quad + \int_{U \cap \{f \ge K \}}\sum_{k}|\eta_{\varepsilon_k}\star(g\xi_k) - g\xi_k|fdx\\
    & \le  K\sum_{k=1}^{\infty} \int_{V_k}|\eta_{\varepsilon_k}\star(g\xi_k) - g\xi_k|dx+4\norm{g}_{\infty}\delta \le K\varepsilon+4\norm{g}_{\infty}\delta .
\end{align*}
Therefore, letting $\varepsilon \rightarrow 0$ and then $\delta \rightarrow 0$, the first claim of the lemma is proven.

To prove the second claim, consider 
\begin{align}\label{e: variation of smoothed g} \nonumber
    \int_{U} g_{\varepsilon}(\text{div}(\phi)f + \phi \cdot \nabla f)dx & = \sum_{k}\int_U \eta_{\varepsilon_k}\star (g\xi_k)[\text{div}(\phi)f + \phi \cdot \nabla f]dx\\
    & = \sum_{k}\int_U g\xi_k[\eta_{\varepsilon_k}\star (\text{div}(\phi)f + \phi \cdot \nabla f)]dx.
\end{align}
Now, we focus on the convolution in equation (\ref{e: variation of smoothed g}).  
\begin{align}\label{e:convolution} \nonumber
& \eta_{\varepsilon_k}\star (\text{div}(\phi)f + \phi \cdot \nabla f)(x)\\ \nonumber
&  = \int \eta_{\varepsilon_k}(y)(\text{div}(\phi)(x-y)f(x-y) + \phi(x-y) \cdot \nabla f(x-y))dy\\ 
&  = (\eta_{\varepsilon_k}\star \text{div}(\phi))(x)f(x) + (\eta_{\varepsilon_k}\star \phi)(x) \cdot\nabla f(x)\\ \nonumber
&  \quad  + \int \eta_{\varepsilon_k}(y)\text{div}(\phi)(x-y)[f(x-y) -f(x)]dy\\ \nonumber
& \ \quad + \int\eta_{\varepsilon_k}(y)\phi(x-y) \cdot [\nabla f(x-y)-\nabla f(x)]dy:=\text{(I)}+\text{(II)}+\text{(III)}.
\end{align}
By integration by parts and estimates on the mollifier $\eta$, we may estimate the middle term in \eqref{e:convolution} as follows.
\begin{align*}
|\text{(II)}|&=\left|\int \eta_{\varepsilon_k}(y)\text{div}(\phi)(x-y)[f(x-y) -f(x)]dy\right|\\
& =\left| \int \nabla \eta_{\varepsilon_k}(y)\cdot \phi(x-y)[f(x-y) -f(x)]dy-\int  \eta_{\varepsilon_k}(y)\phi(x-y)\cdot \nabla f(x-y) dy\right|\\
& \le \int |\nabla \eta_{\varepsilon_k}(y)||f(x-y) -f(x)|dy + \int \eta_{\varepsilon_k}(y)|\nabla f(x-y)|dy\\
& \le c(\eta)\varepsilon_k^{-n-1}\int_{B_{\varepsilon_{k}}(0)} |f(x-y) -f(x)|dy + \varepsilon_k^{-n}\int_{B_{\varepsilon_{k}}(0)} |\nabla f(x-y)|dy\\
& \le c(\eta)\varepsilon_k^{-n}\int_{B_{\varepsilon_{k}}(0)} \int_{0}^1|\nabla f(x-ty)|dtdy + \varepsilon_k^{-n}\int_{B_{\varepsilon_{k}}(x)} |\nabla f(y)|dy.
\end{align*}

Estimating the last term in \eqref{e:convolution} by 
\begin{align*}
   |\text{(III)}|= \left|\int \eta_{\varepsilon_k}(y)\phi(x-y) \cdot [\nabla f(x-y)-\nabla f(x)]dy\right| & \le \varepsilon_k^{-n}\int_{B_{\varepsilon_k}(x)} |\nabla f(y)|dy + |\nabla f(x)|,
\end{align*}
we plug these estimates into \eqref{e: variation of smoothed g} to obtain the following.
\begin{align}\label{four_terms}
    &\left| \int_{U} g_{\varepsilon}(\text{div}(\phi)f + \phi \cdot \nabla f)dx\right| =\left| \sum_{k}\int_U g\xi_k\big[\text{(I)}+\text{(II)}+\text{(III)}\big]dx\right|\notag\\
    &\le \left|\sum_{k}\int_U g\xi_k\Big[(\eta_{\varepsilon_k}\star (\text{div}(\phi))f + (\eta_{\varepsilon_k}\star \phi) \cdot \nabla f)\Big]dx\right|\notag\\
    & \quad +  \sum_{k}\int_U |g|\xi_k\Big[2\varepsilon_k^{-n}\int_{B_{\varepsilon_k}(x)} |\nabla f(y)|dy + |\nabla f(x)|\Big]dx\notag\\
    &\quad +c(\eta)\sum_{k}\int_U |g|\xi_k\Big[\varepsilon_k^{-n}\int_{B_{\varepsilon_{k}}(0)} \int_{0}^1|\nabla f(x-ty)|dtdy\Big]dx\notag\\
    & =\left| \sum_{k}\int_U g[(\text{div}(\xi_k (\eta_{\varepsilon_k}\star\phi))f + \xi_k(\eta_{\varepsilon_k}\star \phi) \cdot \nabla f)]dx - \sum_{k}\int_U g \nabla\xi_k \cdot (\eta_{\varepsilon_k} \star \phi)fdx\right|\\
    & \quad +  \sum_{k}\int_U |g|\xi_k\Big[2\varepsilon_k^{-n}\int_{B_{\varepsilon_k}(x)} |\nabla f(y)|dy + |\nabla f(x)|\Big]dx\notag\\
    &\quad +c(\eta)\sum_{k}\int_U |g|\xi_k\Big[\varepsilon_k^{-n}\int_{B_{\varepsilon_{k}}(0)} \int_{0}^1|\nabla f(x-ty)|dtdy\Big]dx.\notag
\end{align}
Now, we split the first term in equation \eqref{four_terms} into two parts. Because $\xi_k (\eta_{\varepsilon_k} \star \phi) \in C^{1}_c(U)$, $|\xi_k (\eta_{\varepsilon_k} \star \phi)|\le 1$, and $\sum_{k}\chi_{V_k}(x) \le 2$, we have
\begin{align*}
&\left| \sum_{k}\int_U g[(\text{div}(\xi_k (\eta_{\varepsilon_k}\star\phi))f + \xi_k(\eta_{\varepsilon_k}\star \phi) \cdot \nabla f)]dx\right|\\
    & \le \norm{Dg}(U_1) + 2\norm{Dg}(U-U_1) \le \norm{Dg}(U_1) + 2\varepsilon.
\end{align*}
Similarly, because $\sum_k \nabla\xi_k(x) \equiv 0$ we may use \eqref{e:POU 3} to obtain
\begin{align*}
    & \left|\sum_{k}\int_U g \nabla\xi_k \cdot (\eta_{\varepsilon_k} \star \phi)fdx\right| = \left|\sum_{k}\int_U \eta_{\varepsilon_k} \star (gf \nabla\xi_k) \cdot \phi dx - \sum_{k}\int_U g f\nabla\xi_k \cdot \phi dx\right|\\
    & \le  \int_U \sum_{k} |\eta_{\varepsilon_k} \star (gf \nabla\xi_k) - gf\nabla\xi_k| dx \le \varepsilon.
\end{align*}
Finally, using our assumption that $\text{supp}(\eta_{\varepsilon_k} \star \xi_k) \subset V_k$, $\sum_{k}\chi_{V_k}(x) \le 2$, and Fubini's theorem, we estimate the last two terms in equation \eqref{four_terms} as follows. 
\begin{align*}
    & \sum_{k}\int_U |g|\xi_k\Big[2\varepsilon_{k}^{-n}\int_{B_{\varepsilon_{k}}(x)}|\nabla f(y)|dy + |\nabla f(x)|\Big]dx\\
    &\quad+c(\eta)\sum_{k}\int_U |g|\xi_k\Big[\varepsilon_k^{-n}\int_{B_{\varepsilon_{k}}(0)} \int_{0}^1|\nabla f(x-ty)|dtdy\Big]dx\\
    &\leq C(n,\eta)\norm{g}_{\infty}\sum_{k}\int_{V_k}|\nabla f|dx\leq 2C(n,\eta)\norm{g}_{\infty}\int_{U}|\nabla f|dx.
\end{align*}

Putting it all together, we obtain
\begin{align*}
    \left|\int_{U} g_{\varepsilon}(\text{div}(\phi)f + \phi \cdot \nabla f)dx\right| \le \norm{Dg}(U) + 3\varepsilon+ 2C(n,\eta)\norm{g}_{\infty}\int_{U}|\nabla f|dx.
\end{align*}

Letting $\varepsilon \rightarrow 0$ gives the desired result.
\end{proof}

In this paper, we only apply Lemma \ref{l:weak approx by smooth} in the case $U$ is bounded. Though we will not use it, we state the following interesting improvement on Lemma \ref{l:weak approx by smooth} (2).

\begin{lemma}\label{l:weak convergence improvement}
    Under the assumptions of Lemma \ref{l:weak approx by smooth}, suppose additionally that $U$ is bounded. Let $\mathcal{B}$ be the collection of open sets $O \subset \mathbb{R}^n$ such that $g|_{O}$ is constant. Define 
    \begin{align*}
\mathcal{K} = \bigcup\limits_{O \in \mathcal{B}}O.
    \end{align*}
    Then, the conclusion of Lemma \ref{l:weak approx by smooth} (2) may be improved to 
    \begin{align*}
        \limsup_{\varepsilon \rightarrow 0}\norm{Dg_\varepsilon}(U) \le \norm{Dg}(U) + C(n)\norm{g}_{\infty}\norm{\nabla f}_{L^1(U \setminus \mathcal{K})}.
    \end{align*}
    As an immediate corollary, if $g= \chi_E$ for $E$ a set with 
    \begin{align}\label{e:minkoski content of boundary}
        \limsup_{r \rightarrow 0} \mu(B_{r}(\emph{spt}\norm{D\chi_E})) = 0,
    \end{align}
    there exists a sequence of smooth functions $g_{j}$ such that for any (possibly unbounded) open $V \subset \mathbb{R}^n$,
    \begin{align}\label{e:weighted improved approx}
    g_{j} \rightarrow \chi_E \text{ in } L^1_\loc(\mu, V)\quad\mbox{and}\quad
        \limsup_{j \rightarrow \infty}\norm{Dg_j}(\mathbb{R}^n) \le \operatorname{Per}_{\mu}(E; V).
    \end{align}
\end{lemma}

\begin{proof}
First, observe that if $h$ is constant on some open set $O$ then for any $\phi$ with compact support in $O$,
\begin{align}\label{constant_situation}
\int_{\mathbb{R}^n} h(\text{div}(\phi)f + \phi\cdot \nabla f)dx =0.
\end{align}
If $g$ is constant in an open set $O$, then the mollifications $g_{\varepsilon}$ will also be constant in $\{x \in O: \dist(x, \partial O) \ge \sup_{k}\varepsilon_{k}\}$.

Now, we modify the argument from Lemma \ref{l:weak approx by smooth} as follows.  Let $\rho>0$ and choose $\{\varepsilon_k\}_k$ such that additionally $\sup_k \varepsilon_k \le \rho$.  Because $g_{\varepsilon}$ is constant on the components of $U \setminus (B_{\rho}(U \setminus \mathcal{K}))$, by \eqref{constant_situation}, we may reduce to only considering $\phi$ with support compactly contained in $U\cap B_{\rho}(U \setminus \mathcal{K})$ in the definition of $\norm{Dg_{\varepsilon}}$. Making such a choice, we see that \eqref{e: variation of smoothed g} may be replaced by 
    \begin{align*}
        \int_{U\cap B_{\rho}(U \setminus \mathcal{K})} g_{\varepsilon}(\text{div}\phi f+\phi \cdot \nabla f)dx.
    \end{align*}
Carrying this new domain of integration through the rest of the argument in Lemma \ref{l:weak approx by smooth}, we obtain 
\begin{align*}
\limsup_{\varepsilon \rightarrow \infty} \norm{Dg_{\varepsilon}}(U) \le \norm{Dg}(U) + C(n)\norm{g}_{\infty}\norm{\nabla f}_{L^1(U\cap B_{\rho}(U \setminus \mathcal{K}))}.
\end{align*}
We obtain the first claim of this lemma by noting that as $\rho \rightarrow 0$ 
    \begin{align*}
        \int_{U\cap B_{\rho}(U \setminus \mathcal{K})} |\nabla f|dx \rightarrow  \int_{U \setminus \mathcal{K}} |\nabla f|dx
    \end{align*} by monotone convergence.
    
In the special case that $g= \chi_E$ for $E\subset \mathbb{R}^n$, we may assume that $\text{Per}_{\mu}(E;V)< \infty$ otherwise there is nothing to prove.  Let $\text{spt}\norm{D\chi_E}$ be as in the hypotheses.  Let $U_j$ be a sequence of nested bounded, open subsets of $V$ which exhaust $V$. Let $\{O_i\}_i$ be the components of $V\setminus \text{spt}\norm{D\chi_ E}$ and $\mathcal{K}$ as in the statement of the lemma using these $O_i$.  Note that the condition \eqref{e:minkoski content of boundary} implies that 
\begin{align*}
    \int_{U_j \setminus \mathcal{K}} |\nabla f|dx & =\lim_{r \rightarrow 0} \int_{U_j \cap (B_{r}(U_j \setminus \mathcal{K}))} |\nabla f|dx\\
    & =  \lim_{r \rightarrow 0} \int_{U_j \cap (B_{r}(U_j \setminus \mathcal{K})) \cap \{f>0\}} \frac{|\nabla f|}{f}\, d\mu =0.
\end{align*}
Therefore, by the first claim of the lemma, for each $U_j$ we may obtain a smooth function $g_j$ such that 
\begin{align*}
    \norm{g_{j} - g}_{L^1(\mu, U_j)} \le 1/j,\quad   \norm{Dg_j}(U_j) \le \norm{Dg}(U_j) + 1/j.
\end{align*}
Taking the limit $j \rightarrow \infty$ gives the claim.
\end{proof}   

\begin{remark}
    We note that the condition \eqref{e:minkoski content of boundary} does not hold for all sets of finite $\mu$-perimeter.  By ensuring that the radii shrink sufficiently fast, one may obtain a set of finite perimeter from the union of a collection of discs, centered on some dense subset in $B_1(0)\subset \mathbb{R}^2$. See, for example, \cite{Giusti84_Book} Example 1.10. While for this specific example one can obtain the conclusion in \eqref{e:weighted improved approx} by considering finite unions of the discs, it is an interesting question whether or not one can obtain the conclusion in \eqref{e:weighted improved approx} for general sets of finite $\mu$-perimeter.
\end{remark}

\begin{lemma}[Weak Compactness for $BV_{\mu}$]\label{l:weak compactness}
    Let $\mu \in \mathscr{W}(\mathbb{R}^n)$ with distribution function $f$. Let $U \subset \mathbb{R}^n$ be an open set with Lipschitz boundary. For any sequence of functions $h_j \in BV_{\mu}(U)$ such that 
    \begin{align*}
        \sup_j \norm{Dh_j}(U) <\infty, \quad \sup_j \norm{h_j}_{\infty}<\infty,
    \end{align*}
there is a subsequence $\{h_{j'}\}_{j'}$ and a function $h_{\infty} \in BV_{\mu}(U) \cap L^{\infty}(U)$ such that   
\begin{align*}
    h_{j'} \rightarrow h_{\infty} \text{  in  } L^1_{\loc}(\mu,U),\quad  \norm{Dh_{\infty}}(U) \le \liminf_{j\to \infty} \norm{Dh_{j}}(U).
\end{align*} 
In particular, for any sequence of sets of finite $\mu$-perimeter $E_j$ in $U$ such that
$$
\sup _j \operatorname{Per}_\mu\left(E_j ; U\right)<\infty,
$$
there exists a subsequence $\left\{E_{j'}\right\}_{j'}$ and a set of finite $\mu$-perimeter $E_\infty$ such that
$$
\chi_{E_{j'}} \rightarrow \chi_{E_\infty} \text { in } L_{\loc}^1(\mu,U), \quad \operatorname{Per}_\mu(E_\infty ; U) \leq \liminf _{j \rightarrow \infty} \operatorname{Per}_\mu\left(E_{j} ; U\right) .
$$
\end{lemma}

\begin{proof}
    Let $U \subset \mathbb{R}^n$ be given.  For $i \in \mathbb{N}$, we consider the neighborhoods
    \begin{align*}
        U_i = U \cap B_i(0).
    \end{align*}
We claim that for any sequence $\{h_{j}\}$ there is a subsequence $\{h_{j'}\}_{j'}$ and a function $u_i\in BV_{\mu}(U_i) \cap L^{\infty}(U_i)$ such that $h_{j'}\rightarrow u_i$ in $L^1(\mu, U_i)$ and $\norm{Du_i}(U_i) \le \liminf\limits_{j'\to \infty} \norm{Dh_{j'}}(U)$.

Applying Lemma \ref{l:weak approx by smooth} to $\{h_{j}\}_{j}$ in $U_i$, we obtain smooth functions $g_{j} \in BV_{\mu}(U_i)$ such that $\norm{g_{j}}_{\infty} \le (1+2/j)\norm{h_{j}}_{\infty}$ and 
\begin{align}\label{e:smooth approx and strong convergence}
    \begin{cases}
\displaystyle\int_{U_i}|h_{j} - g_{j}|d\mu \le \frac{1}{j}\\
\norm{Dg_{j}}(U_i) \le \norm{Dh_{j}}(U_i) + C(n)\norm{h_{j}}_{\infty}\norm{\nabla f}_{L^1(U_i)}+1.
    \end{cases}
\end{align}
Since we are assuming $\sup_{j}\norm{Dh_{j}}(U_i)<\infty$, $\sup_{j} \norm{h_{j}}_{\infty}<\infty$, and $f \in W^{1,1}_{\loc}(\mathbb{R}^n)$, we have $\sup_{j} \norm{Dg_{j}}(U_i)< \infty$.  That is, there exists some $0 < M <\infty$ such that \begin{align*}
    \sup_{j \in \mathbb{N}} \left(\sup_{\substack{\phi \in C^\infty_c(U_i; \mathbb{R}^n)\\|\phi|\le 1}} \int_{\mathbb{R}^n} (\nabla g_{j} \cdot \phi) fdx\right) < M <\infty.
\end{align*}
For each $j$, we may choose a sequence of test functions $\phi \in C^\infty_c(U_i; \mathbb{R}^n)$ which approximate $\nabla g_{j}/ |\nabla g_{j}|$ in $U_i$ to obtain $\norm{\nabla g_{j}}_{L^1(\mu, U_i)}< M$.

By Leibniz's rule $\nabla (g_{j} f) = (\nabla g_{j}) f + g_{j} (\nabla f)$, and hence the assumptions that $f \in W^{1,1}_{\loc}(\mathbb{R}^n)$ and $\sup_{j}\|g_{j}\|_{\infty} <\infty$ imply that the sequence $\{\nabla (g_{j} f)\}_{j}$ is uniformly bounded in $W^{1,1}(U_i)$. Therefore, by the \textit{unweighted} Rellich-Kondrashov compactness theorem (see, for example, \cite{evansgariepy}, Section 4.6, Theorem 1, Remark), there exists a function $h_{i, \infty} \in L^1(U_i)$ and a subsequence (indexed $j'$) such that
\begin{align*}
    g_{j'}f \rightarrow  h_{i,\infty} \text{ in } L^1(U_i).
\end{align*}

Now, we verify the claim by showing that we can write $h_{i, \infty} = u_{i}f$.  To see this, we note that if $E \subset U_i$ is a set of positive $\mathcal{L}^n$-measure and $f=0$ on $E$, then $\int_E h_{i,\infty} dx = 0$ by the strong convergence.  Therefore, we may define 
\begin{align*}
    u_{i}:= \begin{cases}
        h_{i, \infty}/f & \text{ where } f>0\\
        0 & \text{ where } f=0.
    \end{cases}
\end{align*}

Therefore, $g_{j'}f \rightarrow u_{i}f$ in $L^1(U_i)$ which is the definition of $g_{j'} \rightarrow u_{i}$ in $L^1(\mu, U_i)$.  Recalling \eqref{e:smooth approx and strong convergence} we also have $h_{j'} \rightarrow u_{i}$ in $L^1(\mu, U_i)$.  Thus, we may invoke Lemma \ref{l:lsc} to obtain $\norm{Du_{i}}(U_i) \le \liminf\limits_{j'\to \infty} \norm{Dh_{j'}}(U_i)$.  This proves the claim.

Now, let $\{h_{j'}\}$ be a subsequence of $\{h_j\}$ such that $\lim\limits_{j'\to \infty}\norm{Dh_{j'}}(U) = \liminf\limits_{j\to \infty} \norm{Dh_j}(U)$.  We apply the claim on $U_i$ to obtain a further subsequence $\{h_{k}\}$ and a function $u_i$ on $U_i$.  Applying the claim to this subsequence in $U_{i+1}$ we obtain a further subsequence $\{h_{k'}\}$ and a function $u_{i+1}$ on $U_{i+1}$.  Because $\{h_{k'}\}$ is a subsequence of $\{h_{k}\}$, we must have that $u_i = u_{i+1}$ in $L^1(\mu,U_i) \cap BV_{\mu}(U_i)$. Thus, for each $i \in \mathbb{N}$ we may inductively produce further subsequences $\{h_{k_i}\}$ such that the conclusions of the claim hold and the functions $u_i$ are coherent in the sense that
\begin{align*}
    u_{i+j} = u_i \quad\mathcal{L}^n\text{-a.e. on  } U_i \text{  for all $j \ge 1$}.
\end{align*}
Thus, we may take a diagonalizing sequence of the $\{h_{k_i}\}$ to obtain a subsequence and a function $\lim\limits_{i \rightarrow \infty}u_i := h_{\infty} \in BV_{\mu}(U) \cap L^{\infty}(U)$ such that 
the claims of the lemma holds by Lemma \ref{l:lsc}.

Finally, we apply the previous argument to $h_{j}= \chi_{E_j}$. Then there exists a subsequence such that
$$\mbox{$h_{j'} \rightarrow h_\infty$ in $L^1_\loc(\mu,U)$}.$$
Moreover, there exists a subsequence (again indexed by $j'$) such that $h_{j'}  \rightarrow h_\infty$ $\mu$-a.e. and hence $\mu(\{h_\infty\not=0,1\})=0$, i.e., $h_\infty=\chi_{E_\infty}$ $\mu$-a.e., where $E_\infty=\{h_\infty=1\}$. Hence, by the above claim,
$$\operatorname{Per}_{\mu}(E_\infty; U) \leq \liminf_{j\to \infty} \operatorname{Per}_{\mu}(E_{j}; U).$$
\end{proof}

The next result is a standard result in the theory of Sobolev functions.  We shall use it in the proof of Lemma \ref{l:product structure} as a version of De Giorgi's Structure Theorem for sets of finite $\mu$-perimeter in the special case that the perimeter is the graph of a smooth function with compact support.  We note that the usual proof of De Giorgi's Structure Theorem (see \cite{evansgariepy}, Section 5.7.3, Theorem 2) or the Gauss-Green formula for sets of finite perimeter (see \cite{evansgariepy}, Section 5.8, Theorem 1) rely upon the investigation of blow-ups and the relative isoperimetric inequality.  Since the relative isoperimetric inequality and the weighted Poincar\'{e} inequality may fail for weights $f \in L^1(\mathbb{R}^n) \cap W^{1,1}_{\loc}(\mathbb{R}^n)$ and since Lemma \ref{l:product structure} is sufficient for the application at hand, we do not prove a Gauss-Green formula for general sets of finite $\mu$-perimeter. 

\begin{lemma}[\cite{evansgariepy}; Section 4.3, Theorem 1; Section 5.3, Theorem 2]\label{l:structure theorem}
For any open bounded set $U \subset \mathbb{R}^n$ with Lipschitz boundary, there exists a bounded linear operator 
\begin{align*}
    T:W^{1,1}(U) \rightarrow L^1(\mathcal{H}^{n-1},\partial U)
\end{align*}
such that $Tf$ is defined $\mathcal{H}^{n-1}$-almost everywhere by the formula
\begin{align*}
    Tf(x) = \lim_{r \rightarrow 0} \frac{1}{\omega_nr^n}\int_{B_r(x)} f(y)dy.
\end{align*}
Moreover, for every $\phi \in C^1(\mathbb{R}^n; \mathbb{R}^n)$
\begin{align}
    \int_U \emph{div}(\phi) fdx = - \int_U \phi \cdot \nabla f dx + \int_{\partial U} (\phi \cdot \nu) Tf d\mathcal{H}^{n-1},
\end{align}
where $\nu$ is the unit outward normal to $U$.
\end{lemma}

\begin{definition}
    With a slight abuse of notation, because we have chosen $f \in W^{1,1}(\mathbb{R}^n)$ to be its precise representative, we shall use the notation $f$ for both $f$ and $Tf$, since they agree on $\partial U$ wherever they are defined.
\end{definition}

\subsection{Generalized Ehrhard symmetrization}\label{Generalized Ehrhard symmetrization}

In this section, we generalize the notion of Ehrhard symmetrization from $\gamma_n$ to a large class of finite measures.

\begin{definition}[Generalized Ehrhard Symmetrization]\label{d:gen Ehrh symm}
Let $\mu=f\mathcal{L}^1$ with a non-negative distribution function $f \in L^1(\mathbb{R})$. For a vector $\vec{v} \in \mathbb{S}^{0} = \{-1, 1\}$ we define {\bf generalized $\vec{v}$-Ehrhard symmetrization with respect to $\mu$} to be the set-valued map $S_{\vec{v}}$ described below.  For any $\mathcal{L}^1$-measurable $E \subset \mathbb{R}$, we define
\begin{align*}
S_{\vec{v}}(E):= \{x: \vec{v} \cdot x \ge c \}
\end{align*}
where $c \in \mathbb{R}\cup\{\pm\infty\}$ is chosen to be the \textit{smallest} constant such that $\mu(E) = \mu(S_{\vec{v}}(E))$. Note that this makes $S_{\vec{v}}(E)$ the largest half-line in the $\vec{v}$ direction with $\mu$-measure equal to $\mu(E)$. Sometimes we shall use the notation $E^s_{\vec{v}} = S_{\vec{v}}(E)$.

For $\mu=f\mathcal{L}^n$ with a non-negative distribution function $f \in L^1(\mathbb{R}^n)$ and any vector vector $\vec{v} \in \mathbb{S}^{n-1}$, we define {\bf generalized $\vec{v}$-Ehrhard symmetrization with respect to $\mu$} to be the set-valued function $S_{\vec{v}}$ described below.  For any $\mathcal{L}^n$-measurable set $E \subset \mathbb{R}^n$ and $x \in \vec{v}^{\perp}$, we define $E_x = E \cap \pi_{\vec{v}^{\perp}}^{-1}(x)$ where $\pi_{\vec{v}^{\perp}}$ is the orthogonal projection onto $\vec{v}^{\perp}$.  We define
\begin{align*}
    S_{\vec{v}}(E) := \bigcup_{x \in \vec{v}^{\perp}}(E_x)^s_{\vec{v}},
\end{align*}
where we identify $\pi_{\vec{v}^{\perp}}^{-1}(x)$ with $\mathbb{R}$ and use $\mu_x :=f\mathcal{H}^1\res \pi_{\vec{v}^{\perp}}^{-1}(x)$ to define $(E_x)^s_{\vec{v}}$. For convenience, we shall often use the notation $E^s_{\vec{v}}:= S_{\vec{v}}(E)$.  Because the underlying measure will always be clear from context we omitted it from the notation. 
\end{definition}

\begin{theorem}\label{l:measurability}
    Let $\mu \in \mathscr{W}(\mathbb{R}^n)$ be a measure with distribution function $f$. The generalized $\vec{v}$-Ehrhard symmetrization with respect to $\mu$ defined in Definition \ref{d:gen Ehrh symm} preserves measurability.  That is, for any $\mathcal{L}^n$-measurable set $E \subset \mathbb{R}^n$ and any direction $\vec{v} \in \mathbb{S}^{n-1}$, $S_{\vec{v}}(E)$ is a $\mathcal{L}^n$-measurable set.
\end{theorem}

\begin{remark}
    Theorem \ref{l:measurability} is essential for the Definition \ref{Perimeter symmetrization measures} and the arguments in Section \ref{s:half planes are minimizers}.  However, we defer the proof until Appendix \hyperref[s:app B measurability]{B}.
\end{remark}

\begin{definition}
    Let $\mu \in \mathscr{W}(\mathbb{R}^n)$. For a sequence $\vec{v}_1, ..., \vec{v}_N \in \mathbb{S}^{n-1}$ we shall use the notation
\begin{align}
    E^s_{\vec{v}_1, ..., \vec{v}_N} := S_{\vec{v}_N} \circ ... \circ S_{\vec{v}_1}(E). 
\end{align}
\end{definition}

Now that we have rigorous definitions of $\text{Per}_{\mu}$ and generalized Ehrhard symmetrization with respect to $\mu$ for all $\mu \in \mathscr{W}(\mathbb{R}^n)$, we can give a definition of PS measures with respect to these choices.  

\begin{definition}[Perimeter symmetrization measures]\label{Perimeter symmetrization measures}
Let $n \in \mathbb{N}$ and $\mu \in \mathscr{W}(\mathbb{R}^n)$.  We shall call $\mu$ a {\bf perimeter symmetrization measure} (or a {\bf PS measure}) with respect to generalized Ehrhard symmetrization and $\text{Per}_{\mu}$ if for every $\vec{v} \in \mathbb{S}^{n-1}$ and every set of finite $\mu$-perimeter $E$,
\begin{align}
    \operatorname{Per}_{\mu}(E^s_{\vec{v}}) \le \operatorname{Per}_{\mu}(E).
\end{align}
\end{definition}

\begin{definition}(Half spaces)
    For a vector $\vec{v} \in \mathbb{S}^{n-1}$ and $r \in \mathbb{R} \cup \{\pm\infty\}$, we define
    \begin{align*}
        H(\vec{v}, r) := \{x \in \mathbb{R}^n: x \cdot \vec{v} \ge r\}.
    \end{align*}
For a non-empty $\mathcal{L}^n$-measurable set $E \subset \mathbb{R}^n$, we define 
    \begin{align}\label{e:half space of equal measure}
        H_{\mu}(E, \vec{v}):= H(\vec{v}, r(E, \vec{v})) 
    \end{align}
where $r(E, \vec{v})$ is chosen such that $H(\vec{v}, r(E, \vec{v}))$ is the largest half-space satisfying $\mu(H(\vec{v}, r)) = \mu(E)$.  Note that $H_{\mu}(E, \vec{v})$ only depends upon $\vec{v}, \mu(E)$ and that in particular $H_{\mu}(F, \vec{v}) = H_{\mu}(E, \vec{v})$ for all sets $F \subset \mathbb{R}^n$ such that $\mu(F) = \mu(E).$ 
\end{definition}

\begin{definition}\label{d: symmetric in measure}
Let $f:\mathbb{R}\to \mathbb{R}$. We say that $f$ is {\bf symmetric (around $x$)} if for all $y_1<x<y_2$ with $|y_1-x| = |x-y_2|$, we have $f(y_1) = f(y_2)$. Similarly, we say that a measure $\mu$ on $\mathbb{R}$ is {\bf symmetric (around $x$)} if for all $y_1<x<y_2$ with $|y_1-x| = |x-y_2|$,
$$\mu([y_1, x]) = \mu([x, y_2]).$$
Note that if $\mu =f\mathcal{L}$ with a positive distribution $f\in C(\mathbb{R})\cap L^1(\mathbb{R})$, then $f$ being symmetric (around $x$) implies $\mu$ is also symmetric (around $x$).
\end{definition}

The connection between generalized Ehrhard symmetrization and symmetric measures is well-known. For example, in Bobkov's solution to the characterization problem in $n=1$ for log-concave measures and generalized Ehrhard symmetrization, he proves that in this class, $\nu$ is a PS measure if and only if $\nu$ is symmetric, using the $\varepsilon$-enlargement definition of perimeter (\cite{Bobkov96} Proposition 2.2). 

For the remainder of the paper, when we say that $\mu \in \mathscr{W}(\mathbb{R}^n)$ is a PS measure, we shall mean that it is a PS measure with respect to generalized Ehrhard symmetrization with respect to $\mu$ and weighted $\mu$-perimeter $\text{Per}_{\mu}$.

Since we are not assuming that the distribution function is log-concave, we give the following, related statement.

\begin{lemma}\label{l:Ehrhard implies symmetric}
    Let $\nu$ be a PS measure on $\mathbb{R}$ with distribution $f \in C(\mathbb{R})$ satisfying $f>0$. Then $f$ is symmetric. 
\end{lemma}

\begin{remark}
    We do not specify the notion of perimeter, here, as we shall only need to calculate the perimeter of sets of the form $(-\infty, y]$. The statement holds for any notion of weighted perimeter such that $\text{Per}_{\nu}((-\infty, y]) = f(y)$.
    \end{remark}

\begin{proof}
    For any $y_1,y_2$ such that $\nu((-\infty, y_1]) = \nu([y_2, \infty))$, we have $(-\infty, y_1] = ([y_2, \infty))^s_{-\vec{e_1}}$ and $((-\infty, y_1])^s_{\vec{e_1}} = [y_2, \infty)$.  Therefore, if $\nu$ is a PS measure, then $\operatorname{Per}_{\nu}((-\infty, y_1]) = \operatorname{Per}_{\nu}([y_2,\infty))$ or, equivalently, $f(y_1)=f(y_2)$.

    Now, let $F: \mathbb{R} \rightarrow (0, \nu(\mathbb{R}))$ be defined by $x \mapsto \nu((-\infty, x])$.  Note that $F \in C^1(\mathbb{R})$ and $F^{-1}$ is well-defined since $f>0$.  Moreover, $(F^{-1})'(p) = \frac{1}{f(F^{-1}(p))}$ and the above paragraph implies that for any $p \in (0, \nu(\mathbb{R}))$, $f(F^{-1}(p)) = f(F^{-1}(\nu(\mathbb{R}) -p))$. Thus, $(F^{-1})'(p) = (F^{-1})'(\nu(\mathbb{R})-p)$.  Integrating, we see that there must exist a constant $m$ such that for all $p \in (0, \nu(\mathbb{R}))$,
    \begin{align*}
        F^{-1}(p) + F^{-1}(\nu(\mathbb{R})-p) = 2m.
    \end{align*}
    Therefore, $m$ is the median and $\mu$ is symmetric around $m$. Thus, $f$ is symmetric.
\end{proof}

Now we turn our attention to the PS measure on $\mathbb{R}^n$. Our next lemma shows that if a PS measure possesses a product structure, indicating that it is a product measure $\mu_1 \times \cdots \times \mu_n$, then each component measure $\mu_i$ must be symmetric.

\begin{lemma}\label{l:Ehrhard implies symmetric_2}
    Let $\mu=\mu_1\times \cdots \times \mu_n \in \mathscr{W}(\mathbb{R}^n)$ be a PS measure and let $\mu_i =f_i\mathcal{L}^1$ with distribution $f_i \in C(\mathbb{R})$ satisfying $f_i>0$ for all $i=1,...,n$. Then $f_i$ is symmetric for all $i=1,...,n$. 
\end{lemma}
\begin{proof}
Let $f(x_1,\cdots, x_n)=f_1(x_1)\cdots f_n(x_n)$. The product measure $\mu$ can be expressed as
$$\mu=f\mathcal{L}^n\in \mathscr{W}(\mathbb{R}^n)$$
with $f\in C(\mathbb{R}^n)$ and $f>0$.\\ 
For any $y_1,y_2$ such that $\mu\left((-\infty, y_1]\times \mathbb{R}^{n-1}\right) = \mu\left([y_2, \infty)\times \mathbb{R}^{n-1}\right)$, we have
$$(-\infty, y_1]\times\mathbb{R}^{n-1} = \left([y_2, \infty)\times \mathbb{R}^{n-1}\right)^s_{-e_1},\quad \left((-\infty, y_1] \times \mathbb{R}^{n-1}\right)^s_{e_1}= [y_2, \infty)\times \mathbb{R}^{n-1},$$  
since $\mu=\mu_1\times \cdots \times \mu_n$. Therefore, if $\mu$ is a PS measure, then 
$$\operatorname{Per}_{\mu}\left([y_2,\infty)\times \mathbb{R}^{n-1}\right)=\operatorname{Per}_{\mu}\left(\left((-\infty, y_1] \times \mathbb{R}^{n-1}\right)^s_{e_1}\right) \leq \operatorname{Per}_{\mu}\left((-\infty, y_1]\times\mathbb{R}^{n-1}\right),$$
and
$$ \operatorname{Per}_{\mu}\left((-\infty, y_1]\times\mathbb{R}^{n-1}\right)= \operatorname{Per}_{\mu}\left(\left([y_2, \infty)\times \mathbb{R}^{n-1}\right)^s_{-e_1}\right) \leq \operatorname{Per}_{\mu}\left([y_2,\infty)\times \mathbb{R}^{n-1}\right).$$ 
That is, $\operatorname{Per}_{\mu}\left((-\infty, y_1]\times\mathbb{R}^{n-1}\right)=\operatorname{Per}_{\mu} \left([y_2,\infty)\times \mathbb{R}^{n-1}\right)$. Notice that Theorem \ref{l:structure theorem} and the continuity of $f$ give \begin{align*}
    \operatorname{Per}_{\mu}\left((-\infty, y_1]\times \mathbb{R}^{n-1}\right) &= \int_{\mathbb{R}^{n-1}}f(y_1,x_2,\cdots,x_n)\, dx_2\cdots dx_n\\
&=\int_{\mathbb{R}^{n-1}}f_1(y_1)f_2(x_2)\cdots f_n(x_n)\, dx_2\cdots dx_n\\
&=f_1(y_1)\mu_2(\mathbb{R})\cdots \mu_{n}(\mathbb{R}).
\end{align*}
Similarly, 
$$\operatorname{Per}{\mu}\left([y_2, \infty)\times \mathbb{R}^{n-1}\right)= f_1(y_2)\mu_2(\mathbb{R})\cdots \mu_{n}(\mathbb{R}).$$ Hence, $f_1(y_1)=f_1(y_2)$. We can apply the rest of the proof in Lemma \ref{l:Ehrhard implies symmetric} to $\mu_1$ to conclude that $f_1$ is symmetric.
\end{proof}

The following lemma tells us that shifting a PS measure results in another PS measure. Additionally, we can shift our PS measure with a product structure so that each component is symmetric around zero.

\begin{lemma}\label{l:Ehrhard implies symmetric_3}
    Let $\mu\in \mathscr{W}(\mathbb{R}^n)$ be a PS measure with distribution $f\in C(\mathbb{R}^n)$ satisfying $f>0$ and let $a\in \mathbb{R}^n$. Then
$$\nu:=\widetilde{f}\mathcal{L}^n$$
is also a PS measure, where $\widetilde{f}(x):=f(x+a)$. Moreover, if $f(x_1,...,x_n)=f_1(x_1)\cdots f_n(x_n)$ with $f_i \in C(\mathbb{R})$ satisfying $f_i>0$, then there exists $m\in \mathbb{R}^n$ such that 
$$\nu:=\widetilde{f}\mathcal{L}^n$$
is a PS measure, where $\widetilde{f}(x)=f(x+m)$, $\widetilde{f}(x_1,...,x_n)=\widetilde{f}_1(x_1)\cdots \widetilde{f}_n(x_n)$, $\widetilde{f}_i>0$, $\widetilde{f}_i\in C(\mathbb{R})$, and $\widetilde{f}_i(-t)=\widetilde{f}_i(t)$ for all $i=1,...,n$.
\end{lemma}
\begin{proof}
Given any finite $\nu$-perimeter set $E$ in $\mathbb{R}^n$ and $\vec{u}\in \SS^{n-1}$,
\begin{align*}
\operatorname{Per}_{\nu}(E)&= \sup \left\{\int_E\Big(\widetilde{f}\div\varphi +\varphi \cdot \nabla \widetilde{f}\Big)dx: \varphi\in C^1_c(\mathbb{R}^n; \mathbb{R}^n), |\varphi| \le 1\right\}\\
&= \sup \left\{\int_{E+a}\Big(f(x)(\div\varphi)(x-a) +\varphi(x-a) \cdot \nabla f(x)\Big)dx: \varphi\in C^1_c(\mathbb{R}^n; \mathbb{R}^n), |\varphi| \le 1\right\}\\
&= \sup \left\{\int_{E+a}\Big(f(\div\widetilde{\varphi} +\widetilde{\varphi} \cdot \nabla f\Big)dx: \widetilde{\varphi}\in C^1_c(\mathbb{R}^n; \mathbb{R}^n), |\widetilde{\varphi}| \le 1\right\}=\operatorname{Per}_{\mu}(E+a).
\end{align*}
Moreover, we have
$$E^s_{\nu,\vec{u}}+a=\left(E+a\right)^s_{\mu,\vec{u}}.$$
Combining them together,
$$\operatorname{Per}_{\nu}(E^s_{\nu,\vec{u}})=\operatorname{Per}_{\mu}\left(E^s_{\nu,\vec{u}}+a\right)=\operatorname{Per}_{\mu}\left(\left(E+a\right)^s_{\mu,\vec{u}}\right)\leq \operatorname{Per}_{\mu}(E+a)=\operatorname{Per}_{\nu}(E)$$
since $\mu$ is a PS measure. Thus, $\nu$ is also a PS measure. In particular, if $f(x_1,...,x_n)=f_1(x_1)\cdots f_n(x_n)$, by Lemma \ref{l:Ehrhard implies symmetric_2}, there exist $m_1,...,m_n$ in $\mathbb{R}$ such that $f_i$ is symmetric around $m_i$ for all $i=1,...,n$. Let $m=(m_1,...,m_n)\in \mathbb{R}^n$ and define $\widetilde{f}(x)=f(x+m)$. That is, $\widetilde{f}(x_1,...,x_n)=\widetilde{f}_1(x_1)\cdots \widetilde{f}_n(x_n)$, where $\widetilde{f}_i(t)=f_i(t+m_i)$. Since $f_i$ is symmetric around $m_i$, we have
$$\widetilde{f}_i(-t)=f_i(m_i-t)=f_i(m_i+t)=\widetilde{f}_i(t)$$
for any $-t<0<t$.
\end{proof}

\section{Half-planes are isoperimetric sets}\label{s:half planes are minimizers}

The main result in this Section is the following lemma.

\begin{lemma}\label{l:H are minimal}
    Let $n \ge 2$, and let $\mu \in \mathscr{W}(\mathbb{R}^n)$ be a PS measure. For every $\vec{v} \in \mathbb{S}^{n-1}$ and every measurable set $E \subset \mathbb{R}^n$,
    \begin{align*}
    \operatorname{Per}_{\mu}(H_{\mu}(E,\vec{v})) \le \operatorname{Per}_{\mu}(E).
    \end{align*}    
\end{lemma}

The main idea of the proof is that for any $\mathcal{L}^n$-measurable set $E \subset \mathbb{R}^n$ if $\mu$ is a PS measure then generalized Ehrhard symmetrization allows us to move mass from $E \setminus H_{\mu}(E, \vec{v})$ to ``fill in" $H_{\mu}(E, \vec{v}) \setminus E$ while not increasing the perimeter (see Figure \ref{fig:moving mass}). As noted in the Introduction, this type of argument may be found for Steiner symmetrization in Euclidean space (\cite{BuragoZalgaller_GeometricInequalitiesBOOK} Lemma 9.4.3), but differs from Ehrhard's original proof of Lemma \ref{l:H are minimal} for $\mu = \gamma_n$ in \cite{Ehrhard83} Proposition 1.5 (again, see Lifschits's book \cite{lifshits1995gaussian}, Chapter 11 for a detailed proof in English).  The core of Ehrhard's proof is to produce a sequence of symmetrizations and show that $C_N \subset E^s_{\vec{v_1}\vec{v_2}...\vec{v_N}}$ for a sequence $\{C_i\}_{i=1}^{\infty}$ of cones which widen to the half-space.  This proof relies essentially upon the measure $\mu$ having a product structure for any orthonormal frame. While this is obviously true for the isotropic Gaussian, we cannot assume this for general $\mu \in \mathscr{W}$.

We begin by studying the one-dimensional case. 

\begin{lemma}[One-Dimensional Behavior]\label{l:0}
    Let $\mu =f\mathcal{L}^1$ for some $f \in L^1(\mathbb{R})$ such that $f\ge 0$. Suppose $E \subset \mathbb{R}$ is a $\mathcal{L}^1$-measurable set and the interval $(a, b)$ satisfies 
    \begin{align*}
        \mu(E^c \cap (a, b)) = \alpha,\quad \mu(E\cap [b,\infty))=\beta.
    \end{align*}
    Then $\mu \left(E^s_{-\vec{e_1}} \cap [b,\infty)\right) \le \max\{ \beta - \alpha, 0\}$. Equivalently, 
$$\mu \left(E^s_{-\vec{e_1}} \setminus (-\infty,b]\right) \le \mu \left(E \setminus (-\infty,b]\right) - \min \{  \beta, \alpha\}.$$  
In particular,
\begin{enumerate}
\item[(1)] $\mu\left(E_{-e_1}^s \setminus(-\infty, b]\right) \leq \mu(E \setminus (-\infty, b])$.
\item[(2)] if there exists $x>y$, and $r>0$ such that $(x-r, x+r) \cap(y-r, y+r)=\emptyset$ and
$$
\mu(E \cap(x-r, x+r))>0, \quad \mu\left(E^c \cap(y-r, y+r)\right)>0,
$$
then $\mu\left(E_{-e_1}^s \setminus(-\infty, b]\right)<\mu(E \setminus(-\infty, b])$ for any $y+r<b<x-r$.
\end{enumerate}
\end{lemma}
\begin{proof}
Without loss of generality, we may assume that $E^s_{-\vec{e_1}} \cap [b,\infty)\not=\emptyset$. Since $E^s_{-\vec{e_1}}$ is an interval of the form $(-\infty,c]$ for some $c$, we have  $E^s_{-\vec{e_1}}\cup [b,\infty)=\mathbb{R}$. Notice that
$$\left(E^c\cap (a,b)\right)\cap \left(E\cup [b,\infty)\right)=\emptyset.$$
Then
\begin{align*}
\alpha+\mu \left(E\cup [b,\infty)\right)&=\mu\left(E^c\cap (a,b)\right)+\mu \left(E\cup [b,\infty)\right)\\
&=\mu\left((E^c\cap (a,b))\cup (E\cup [b,\infty)) \right)\leq \mu\left(\mathbb{R}\right)
\end{align*}
and hence
\begin{align*}
\mu \left(E^s_{-\vec{e_1}} \cap [b,\infty)\right)&=\mu \left(E^s_{-\vec{e_1}}\right)+\mu \left( [b,\infty)\right)-\mu \left(E^s_{-\vec{e_1}} \cup [b,\infty)\right)\\
&=\mu \left(E\right)+\mu \left( [b,\infty)\right)-\mu \left(\mathbb{R}\right)\\
&=\mu\left(E\cap [b,\infty)\right)+\mu \left(E\cup [b,\infty)\right)-\mu \left(\mathbb{R}\right)\\
&=\beta-\left(\mu\left(\mathbb{R}\right)-\mu \left(E\cup [b,\infty)\right)\right)\leq \beta-\alpha.
\end{align*}
Observe that $\mu(E \setminus(-\infty, b]))=\mu(E \cap[b, \infty))=\beta$ since $\mu(\{b\})=0$. This gives us,
$$
\begin{aligned}
\mu\left(E_{-\vec{e_1}}^s \setminus(-\infty, b]\right)=\mu\left(E_{-\vec{e_1}}^s \cap[b, \infty)\right) &\leq \max \{\beta-\alpha, 0\}  =\beta-\min \{\beta, \alpha\} \\
& =\mu(E \setminus(-\infty, b])-\min \{\beta, \alpha\}
\end{aligned}
$$
\end{proof}

Lemma \ref{l:0} has the following, immediate corollary in higher dimensions.

\begin{lemma}\label{l:1}
Let $n \geq 2$, $\mu=f \mathcal{L}^n \in \mathscr{W}\left(\mathbb{R}^n\right)$ be a PS measure, $E \subset \mathbb{R}^n$ be a $\mathcal{L}^n$-measurable set, and $\vec{v} \in \mathbb{S}^{n-1}$. Then for every $\vec{\eta} \in \mathbb{S}^{n-1}$ with $\vec{v} \cdot \vec{\eta}>0$,
\begin{align}
    \mu(E^s_{\vec{\eta}} \setminus H_{\mu}(E, \vec{v})) \le \mu(E \setminus H_{\mu}(E, \vec{v})).
\end{align}
\end{lemma}

\begin{proof}
    This follows from Fubini's theorem and Lemma \ref{l:0} (1) applied to each of the fibers $\pi^{-1}_{\vec{\eta}^\perp}(z)$ for $z \in \langle \vec{\eta}\rangle^\perp$ with measure $\mu_z=f\mathcal{H}^{1}\res \pi_{\vec{\eta}^{\perp}}^{-1}(z)$, where $\pi_{\vec{\eta}^{\perp}}$ is the orthogonal projection onto $\langle\vec{\eta}\rangle^{\perp}$. The following figure illustrates how we use Lemma \ref{l:0}  to move mass through fibers (1) and (2).
    
\begin{figure}[!ht]
    \centering
    \includegraphics[width=.74\textwidth]{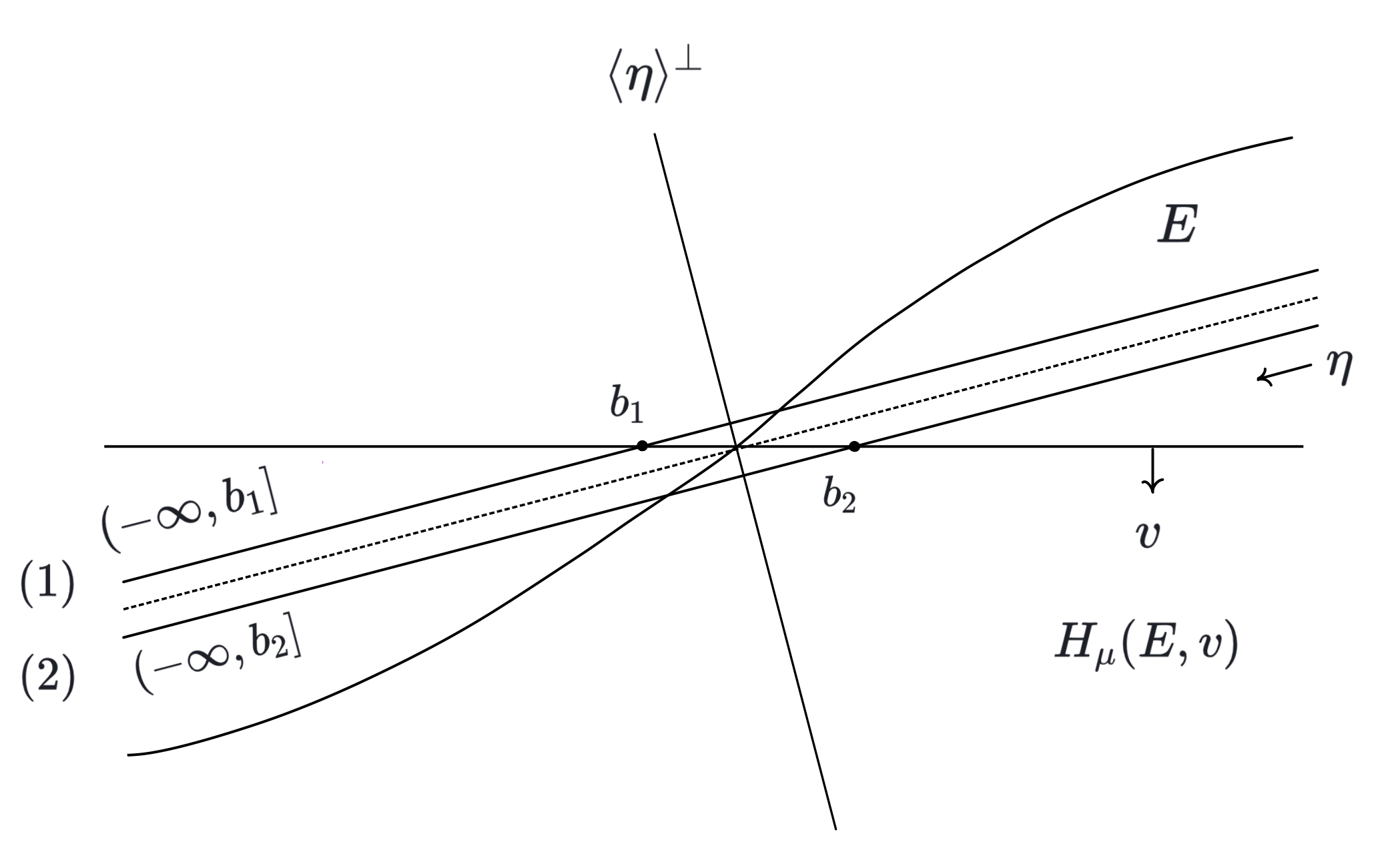}
    \caption{Moving mass from $E \setminus H_{\mu}(E, \vec{v})$ to ``fill in" $H_{\mu}(E, \vec{v}) \setminus E$.}
    \label{fig:moving mass}
\end{figure}
\end{proof}

However, Lemma \ref{l:1} is not strong enough to prove Lemma \ref{l:H are minimal}. Instead, we need the following improvement.

\begin{lemma}\label{l:measure reduction}
Let $n \geq 2$, and let $\mu=f \mathcal{L}^n \in \mathscr{W}\left(\mathbb{R}^n\right)$ be a PS measure, $E \subset \mathbb{R}^n$ be a $\mathcal{L}^n$-measurable set, $\vec{v} \in \mathbb{S}^{n-1}$, and assume that
$$
\mu\left(E \Delta H_\mu(E, \vec{v})\right)>0.
$$
There exists a vector $\vec{\eta} \in \mathbb{S}^{n-1}$ such that $\vec{v} \cdot \vec{\eta}>0$ and 
\begin{align}
    \mu(E^s_{\vec{\eta}} \setminus H_{\mu}(E, \vec{v})) <  \mu(E \setminus H_{\mu}(E, \vec{v})). 
\end{align}
In particular, $\mu(E^s_{\vec{\eta}} \Delta H_{\mu}(E, \vec{v})) <  \mu(E \Delta H_{\mu}(E, \vec{v})). $
\end{lemma}

\begin{proof}
By the assumption that $\mu(E\Delta H_{\mu}(E, \vec{v}))>0$, we may find an $\varepsilon>0$ such that 
\begin{align*}
    \mathcal{H}^{n}(\{f \ge \varepsilon\} \cap \left(E \setminus H_{\mu}(E, \vec{v}))\right) & > 0\\
    \mathcal{H}^n(\{f \ge \varepsilon\} \cap \left(H_{\mu}(E, \vec{v}) \setminus E\right)) & > 0.
\end{align*}
Therefore, taking $x, y$ such that 
\begin{enumerate}
    \item[(a)] $x$ is a point of $\mathcal{L}^n$-density 1 for $\{f \ge \varepsilon\}\cap (E \setminus H_{\mu}(E, \vec{v}))$.
    \item[(b)] $y$ is a point of $\mathcal{L}^n$-density 1 for $\{f \ge \varepsilon\}\cap (H_{\mu}(E, \vec{v})\setminus E)$.
    \item[(c)] $\vec{v}\cdot (y-x)>0$.
\end{enumerate}
Let $\vec{\eta} = \frac{y-x}{|y-x|}$ and let $\vec{\eta}^\perp_x := x+ \vec{\eta}^\perp$. For $z \in \vec{\eta}^{\perp}_x$, we denote $\mu_z=f\mathcal{H}^{1}\res \pi_{\vec{\eta}_x^{\perp}}^{-1}(z)$, where $\pi_{\vec{\eta}_x^{\perp}}$ is the orthogonal projection onto $\vec{\eta}^{\perp}_x$.  Note that there exists a radius $r_0>0$ such that 
\begin{align*}
    \frac{\mathcal{H}^{n-1}(\{z \in \vec{\eta}^\perp_x \cap B_{r_0}(x): \mu_z(E \cap B_{r_0}(x))>0\})}{\mathcal{H}^{n-1}(\vec{\eta}^\perp_x \cap B_{r_0}(x))} & > 3/4\\
    \frac{\mathcal{H}^{n-1}(\{z \in \vec{\eta}^\perp_x \cap B_{r_0}(x): \mu_z(E^c \cap B_{r_0}(y))>0\})}{\mathcal{H}^{n-1}(\vec{\eta}^\perp_x \cap B_{r_0}(x))} & > 3/4.
\end{align*}
For this $r_0>0$ the collection of $z \in \vec{\eta}^\perp_x \cap B_{r_0}(x)$ such that both $\mu_z(E \cap B_{r_0}(x))>0$ and $\mu_z(E^c \cap B_{r_0}(y))>0$ must have positive $\mathcal{H}^{n-1}$-measure. Therefore, applying Lemma \ref{l:0} (2) to each of the fibers inside the tube based on the above set from $B_{r_0}(x)$ to $B_{r_0}(y)$, and Lemma \ref{l:0} (1) to each of the fibers outside that tube, yields the claim of the lemma.
\end{proof}

\begin{lemma}\label{l:filling all the holes}
Let $n \geq 2, \mu=f \mathcal{L}^n \in \mathscr{W}\left(\mathbb{R}^n\right)$ be a PS measure, $\vec{v} \in \mathbb{S}^{n-1}$, and assume that $E$ is a $\mathcal{L}^n$-measurable set with finite $\mu$-perimeter. Let $\mathcal{S}(E)$ be the collection of all sets $F \subset \mathbb{R}^n$ such that 
\begin{align*}
    F = E^s_{\vec{\eta}_1, ..., \vec{\eta}_k}
\end{align*}
for some finite sequence of vectors $\{\vec{\eta}_i \}_{i=1}^k$ such that $\vec{v} \cdot \vec{\eta}_i >0$ for all $i=1,...,k$.  Then
    \begin{align*}
        \inf_{F \in \mathcal{S}(E)} \mu(F \Delta H_{\mu}(E, \vec{v})) =0.
    \end{align*}
\end{lemma}

\begin{proof}
We argue by contradiction.  Suppose that $\inf\limits_{F \in \mathcal{S}(E)} \mu(F \Delta H_{\mu}(E, \vec{v}))>0.
$ Then, for a minimizing sequence $F_j$, we may invoke Lemma \ref{l:weak compactness} to obtain a set $F_{\infty} \subset \mathbb{R}^n$ such that there exists a subsequence (still denoted as $F_j$)
\begin{align*}
    \chi_{F_{j}} \rightarrow \chi_{F_{\infty}} \ \text{ in } L^1_{\loc}(\mu, \mathbb{R}^n),\quad \operatorname{Per}_{\mu}(F_{\infty}) \le \liminf_{j\to \infty}\operatorname{Per}_\mu(F_j),
\end{align*}
where we have used $\operatorname{Per}_\mu\left(F_j\right) \leq \operatorname{Per}_\mu(E)$ since $\mu$ is a PS measure and hence
$$
\sup _j \operatorname{Per}_\mu\left(F_j\right) \leq \operatorname{Per}_\mu(E)<\infty.
$$
Moreover, by the strong convergence in $L^1_{\loc}(\mu, \mathbb{R}^n)$, it also follows that 
\begin{align*}
    \mu(F_{\infty} \Delta H_{\mu}(E, \vec{v}))>0.
\end{align*}
Therefore, we may apply Lemma \ref{l:measure reduction} to $\mu, f, F_\infty, \vec{v}$ to obtain a vector $\vec{\eta} \in \mathbb{S}^{n-1}$ with $\vec{v}\cdot \vec{\eta}>0$ and 
\begin{align*}
    \mu((F_{\infty})^s_{\vec{\eta}} \Delta H_{\mu}(E, \vec{v}))< \inf_{F \in \mathcal{S}(E)} \mu(F \Delta H_{\mu}(E, \vec{v})).
\end{align*}
Therefore, the proof of the lemma will be accomplished if we can show that for sufficiently large $j$ we have $\mu((F_{j})^s_{\vec{\eta}} \Delta H_{\mu}(E, \vec{v})) \rightarrow \mu((F_{\infty})^s_{\vec{\eta}} \Delta H_{\mu}(E, \vec{v}))$.  However, this is immediate from the strong $L^1_{\loc}(\mu, \mathbb{R}^n)$ convergence. That is, for the fixed $\vec{\eta}$ above, $\chi_{F_{j}} \rightarrow \chi_{F_{\infty}}$ in  $L^1_{\loc}(\mu, \mathbb{R}^n)$ implies by Fubini's theorem that the marginal function $m_{F_{j}}: \vec{\eta}^{\perp} \rightarrow \mathbb{R}$ defined by 
\begin{align*}
    m_{F_{j}}(a) & = \int_{\pi_{\vec{\eta}}^{-1}(a) \cap F_{j}} f d\mathcal{H}^1
\end{align*}
satisfies $m_{F_{j}} \rightarrow m_{F_{\infty}}$ in $L^1_{\loc}(\vec{\eta}^{\perp})$.  Since the symmetrizations $(F_{j})^s_{\vec{\eta}}, (F_{\infty})^s_{\vec{\eta}}$ are uniquely determined by the functions $m_{F_{j}}, m_{F_{\infty}}$ and 
\begin{align*}
    \mu((F_{j})^s_{\vec{\eta}} \Delta H_{\mu}(E, \vec{v})) & = \int_{\vec{\eta}^\perp} |m_{F_{j}}(a) - m_{H_{\mu}(E, \vec{v})}(a)| d\mathcal{H}^{n-1}(a)\\
    & \le \int_{\vec{\eta}^\perp} |m_{F_{j}}(a) - m_{F_{\infty}}(a)| + |m_{F_{\infty}}(a) -  m_{H_{\mu}(E, \vec{v})}(a)| d\mathcal{H}^{n-1}(a)\\
    & \le \int_{\vec{\eta}^\perp} |m_{F_{j}}(a) - m_{F_{\infty}}(a)|d\mathcal{H}^{n-1}(a) + \mu((F_{\infty})^s_{\vec{\eta}} \Delta H_{\mu}(E, \vec{v})),
\end{align*}
letting $j \rightarrow \infty$, we obtain $\lim\limits_{j \rightarrow \infty} \mu((F_{j})^s_{\vec{\eta}} \Delta H_{\mu}(E, \vec{v})) \le \mu((F_{\infty})^s_{\vec{\eta}} \Delta H_{\mu}(E, \vec{v}))$. Hence, for sufficiently large $j$, there is a $F_{j}$ such that 
\begin{align*}
    \mu((F_{j})^s_{\vec{\eta}} \Delta H_{\mu}(E, \vec{v}))< \inf_{F \in \mathcal{S}(E)} \mu(F \Delta H_{\mu}(E, \vec{v})).
\end{align*}
But, since $(F_{j})^s_{\vec{\eta}} = (E^s_{\vec{\eta}_1,...,\vec{\eta}_k})^s_{\vec{\eta}} \in \mathcal{S}(E)$, we obtain a contradiction.  Therefore, the lemma holds.
\end{proof}

\subsection{Proof of Lemma \ref{l:H are minimal}}
Let $E, \mu, \vec{v}$ be as in the hypotheses of Lemma \ref{l:H are minimal}. If $\operatorname{Per}_\mu(E) = \infty$, there is nothing to prove.  Therefore, we may assume that $E$ is a set of finite $\mu$-perimeter. By Lemma \ref{l:filling all the holes}, we may take a minimizing sequence $F_{j} \in \mathcal{S}(E)$ such that 
\begin{align*}
    \lim_{j \rightarrow \infty} \mu(F_{j} \Delta H_{\mu}(E, \vec{v}))=0.
\end{align*}
Since $\sup_j \operatorname{Per}_{\mu}(F_j) \le \operatorname{Per}_{\mu}(E)<\infty$, we may apply Lemma \ref{l:weak compactness} to extract a subsequence (still denoted as $F_j$) such that $\chi_{F_{j}} \rightarrow \chi_{F_{\infty}}$ in $L^1_{\loc}(\mu, \mathbb{R}^n)$ for some $\mathcal{L}^n$-measurable set $F_{\infty}$ with 
\begin{align*}
    \operatorname{Per}_{\mu}(F_{\infty}) \le \operatorname{Per}_\mu(E).
\end{align*}
Moreover, by $\chi_{F_{j}} \rightarrow \chi_{F_{\infty}}$ in $L^1_{\loc}(\mu, \mathbb{R}^n)$, we also have $\mu(F_{\infty} \Delta H_{\mu}(E, \vec{v})) = 0$. This implies $\operatorname{Per}_{\mu}(F_{\infty}) = \operatorname{Per}_{\mu}(H_{\mu}(E, \vec{v}))$, which concludes the proof. \qed

\vspace{.5cm}
Lemma \ref{l:H are minimal} has the following trivial consequence.

\begin{lemma}\label{l:all half-spaces}
    Let $n \geq 2$, and let $\mu=f \mathcal{L}^n \in \mathscr{W}\left(\mathbb{R}^n\right)$ be a PS measure. For any $\vec{v}, \vec{\eta}$ in $\mathbb{S}^{n-1}$ and any measurable set $E \subset \mathbb{R}^n$,
    \begin{align*}
        \emph{Per}_{\mu}(H_{\mu}(E, \vec{\eta})) = \emph{Per}_{\mu}(H_{\mu}(E, \vec{v}))
    \end{align*}
    and the half-spaces $\{H_{\mu}(E, \vec{\eta})\}_{\vec{\eta} \in \mathbb{S}^{n-1}}$ are minimal for $\mu$-perimeter among all measurable sets $F$ with $\mu(F) = \mu(E)$. In particular, for any $\vec{u}\in \SS^{n-1}$ and $c \in \mathbb{R}$, the half-space $H(u, c)$ is minimal for $\mu$-perimeter among all measurable sets $F$ with $\mu(F)=\mu(H(u, c))$.
\end{lemma}

\begin{proof}
    Applying Lemma \ref{l:H are minimal} to $H_{\mu}(E, \vec{v})$ and $\vec{\eta}$ in place of $\vec{v}$ shows that 
    \begin{align*}
        \operatorname{Per}_{\mu}(H_{\mu}(E, \vec{\eta})) \le \operatorname{Per}_{\mu}(H_{\mu}(E, \vec{v})).
    \end{align*}
    Applying it again gives the reverse inequality.   
\end{proof}

\section{The density \texorpdfstring{$f$}{f} enjoys a product structure}\label{s:product structure}

The main result in this section is to show that if $\mu \in \mathscr{W}(\mathbb{R}^n)$ with distribution function $f$ is a PS measure, then the function $f$ enjoys a product structure subordinate to all orthogonal frames $\{\vec{e_i}\}_i^n$ in $\mathbb{R}^n$. This is made rigorous in Corollary \ref{c:product structure in all directions}.

The key idea is to use variational techniques employed by \cite{Rosales14} in the study of isoperimetric sets for perturbations of log-concave measures and \cite{BrockChiacchioMercaldo12}\cite{BrockChiacchio16} in the study of the structure of measures with foliations by isoperimetric sets. This technique is carried out in Lemma \ref{l:product structure}, and relies essentially upon Lemma \ref{l:structure theorem}. 

\begin{lemma}\label{l:product structure}
    Let $\mu \in \mathscr{W}(\mathbb{R}^n)$ with distribution function $f$. Write $x \in \mathbb{R}^n$ as $x = (x', x_n) \in \mathbb{R}^{n-1}\times \mathbb{R}$.  If the sets $H(\vec{e_n}, c) =\{x \in \mathbb{R}^n: x_n \ge c\}$ are minimal for $\mu$-perimeter among all measurable sets $F$ with $\mu(F)=\mu(H(\vec{e_n}, c))$ for all $c \in \mathbb{R}$, then the following statements hold.
    \begin{enumerate}
        \item There exists a disjoint collection of open intervals $\{I_j\}_{j\in \mathbb{N}}$ such that for $\mathcal{H}^{n-1}$-a.e. $x' \in \mathbb{R}^{n-1}$, the positivity set $\{f(x', \cdot)>0\} = \bigcup_{j \in A}I_j$ for some index set $A \subset \mathbb{N}$ depending on $x'$.
        \item For each $I_k \in \{I_j \}_{j\in \mathbb{N}}$, there exist functions $A_k: \mathbb{R}^{n-1}\rightarrow \mathbb{R}$, $B_k:I_k\rightarrow \mathbb{R}$ such that for $\mathcal{H}^{n}$-a.e. $(x', x_n) \in \mathbb{R}^{n}$,
    \begin{align}\label{e: product structure x' x n}
        f(x',x_n) = A_k(x')B_k(x_n).
    \end{align}
    \item For every $I_k \in \{I_j \}_{j\in \mathbb{N}}$, we have $B_k \in W^{1, 1}_{loc}(I_k) \subset C(I_k)$ and $B_k>0$ on $I_k$.
    \end{enumerate}
    \end{lemma}

\begin{proof}
{\bf Step 1}: Given a function $h \in C^1_c(\mathbb{R}^{n-1})$, we shall construct a corresponding family of measure-preserving perturbations of $H(\vec{e_n}, c)$ as follows. Let $s_h: \mathbb{R}^{n-1}\times (-\varepsilon_0, \varepsilon_0) \rightarrow \mathbb{R}$ be a ``shift" function defined in the following manner.  Given $h\in C^1_c(\mathbb{R}^{n-1})$, we let $s_h \in C^\infty_c(\mathbb{R}^{n-1})$ be a smooth function such that
    \begin{align*}
        s_h(x') = 1 \text{ in a neighborhood of $\text{spt}(h)$,}\\
        s_h(x') \in [0,1] \text{ for all }x' \in \mathbb{R}^{n-1}.
    \end{align*}
    Now, define the real number $c(\varepsilon)$ to be the smallest in absolute value such that the set
    \begin{align*}
G_\varepsilon:= \{x\in \mathbb{R}^n: x_n \ge c+ \varepsilon h(x') + c(\varepsilon)s_h(x') \}\end{align*}
satisfies 
\begin{align}\label{e:measure preserving}
\mu(H(\vec{e_n}, c)) = \mu(G_{\varepsilon}).
\end{align}
That is, for each $\varepsilon>0$, the function $\varepsilon h(x') + c(\varepsilon)s_h(x')$ defines a measure-preserving perturbations of $H(\vec{e_n}, c)$.

We define $S_h(x', \varepsilon)= c(\varepsilon) s_h(x')$. Note that the equality of measure condition \eqref{e:measure preserving} forces $c(0)=0$.  Moreover, if we employ the notation 
$$g(x',\varepsilon):= c+ \varepsilon h(x') + S_h(x',\varepsilon),$$ 
then the condition \eqref{e:measure preserving} can be re-written as
\begin{align}\label{e:meas preserving formula}
\mu(H(\vec{e_n}, c))=\int_{\mathbb{R}^{n-1}}\int_{g(x', \varepsilon)}^{\infty}f(x', x_n)dx_n dx' .
\end{align}
Without loss of generality, we may assume that there is a choice of $s_h$ such that 
$$\int_{\mathbb{R}^{n-1}}s_h(x')f(x', c)dx'\not=0.$$ 
Otherwise, $f(\cdot, c) = 0$ $\mathcal{H}^{n}$-a.e. on $\mathbb{R}^{n-1}$ and \eqref{e: product structure x' x n} holds trivially with $B(x_n) = 0$.
Now, we consider the limit of the difference quotients.
\begin{align*}
0 & = \lim_{\varepsilon \rightarrow 0}\frac{\mu(G_{\varepsilon}) - \mu(H(\vec{e_n}, c))}{\varepsilon} \\
& = \lim_{\varepsilon \rightarrow 0}\frac{\int_{\mathbb{R}^{n-1}}\Big[\int_{g(x', \varepsilon)}^{\infty}f(x', x_n) - \int_{c}^{\infty}f(x', x_n)dx_n \Big] dx'}{\varepsilon}\\
& = \lim_{\varepsilon \rightarrow 0}\frac{\int_{\mathbb{R}^{n-1}}\int_{g(x', \varepsilon)}^{c}f(x', x_n)dx_ndx'}{\varepsilon}.
\end{align*}
Notice that for $\mathcal{H}^{n-1}$-a.e. $x' \in \mathbb{R}^{n-1}$, the function $x_n\mapsto f(x', x_n)$ is in $W^{1,1}_{\loc}(\mathbb{R})$ and hence continuous (see, for example, \cite{brezis} Theorem 8.2). We may use the mean value theorem and find a $x_{n, \varepsilon} \in [\min\{c, g(x', \varepsilon)\}, \max\{c, g(x', \varepsilon)\}]$ such that
\begin{align*}
0=\lim_{\varepsilon \rightarrow 0}\frac{\int_{\mathbb{R}^{n-1}}\int_c^{g(x', \varepsilon)}f(x', x_n)dx_ndx'}{\varepsilon}
& = \lim_{\varepsilon \rightarrow 0}\frac{\int_{\mathbb{R}^{n-1}}(\varepsilon h(x') + S_h(x',\varepsilon))f(x', x_{n, \varepsilon})dx'}{\varepsilon}\\
& = \int_{\mathbb{R}^{n-1}}  f(x', c)h(x') dx' \\
&\quad+ \lim_{\varepsilon \rightarrow 0} \frac{c(\varepsilon) -c(0)}{\varepsilon}\int_{\mathbb{R}^{n-1}}s_h(x')f(x', c)dx',
\end{align*}
where we have used $c(0)=0$. In particular, if $\int_{\mathbb{R}^{n-1}}f(x', c)h(x')dx'=0$, then 
\begin{align}\label{s_prime}
c'(0)=\lim_{\varepsilon \rightarrow 0} \frac{c(\varepsilon) -c(0)}{\varepsilon}=0.
\end{align}
For any $k \in \mathbb{N}$ and $x \in \mathbb{R}^k$, let $Q_r^k(x) \subset \mathbb{R}^k$ be the axis-parallel cube centered at $x$ with side length $r>0$.  By Lemma \ref{l:structure theorem} for any sufficiently large open cube $Q= Q^n_R(0)$ such that $\text{spt}(h) \subset \text{spt}(s_h)\subset B^{n-1}_r(0)$ and $r + |c| < R/2$
\begin{align}\label{e:minimizer}
    \operatorname{Per}_{\mu}(G_\varepsilon; Q) = \int_{\mathbb{R}^{n-1} \cap Q^{n-1}_R(0)}f(x', g(x', \varepsilon))\sqrt{1+ |\nabla_{x'}(g(x', \varepsilon))|^2}\, dx'.
\end{align}
And, by our assumption that $H(\vec{e_n}, c)$ is minimal for $\mu$-perimeter in its mass class,
\begin{align}\label{e:first variation}
     \lim_{\varepsilon\to 0}\frac{\operatorname{Per}_{\mu}(G_{\varepsilon}; Q)-\operatorname{Per}_{\mu}(H(\vec{e_n}, c); Q)}{\varepsilon}= \lim_{\varepsilon\to 0}\frac{\operatorname{Per}_{\mu}(G_{\varepsilon})-\operatorname{Per}_{\mu}(H(\vec{e_n}, c))}{\varepsilon}= 0
\end{align}
since both $h$ and $s_h$ have compact support. By the definition of $s_h$, we have 
$$\nabla h(x') \cdot \nabla s_h(x') = 0$$ 
pointwise. Therefore,
\begin{align*}
0&=\lim_{\varepsilon\to 0}\frac{\operatorname{Per}_{\mu}(G_{\varepsilon}; Q)-\operatorname{Per}_{\mu}(H(\vec{e_n}, c);Q)}{\varepsilon}\\
&=\lim_{\varepsilon\to 0}\frac{1}{\varepsilon}\left(\int_{\mathbb{R}^{n-1} \cap Q^{n-1}_R(0)}f(x', g(x', \varepsilon))\sqrt{1+ |\nabla_{x'}g(x',\varepsilon)|^2}dx'-\int_{\mathbb{R}^{n-1}\cap Q^{n-1}_R(0)}f(x', c)dx'\right)\\
&=\lim_{\varepsilon\to 0}\frac{1}{\varepsilon}\int_{\mathbb{R}^{n-1} \cap Q^{n-1}_R(0)}\Big(f(x', c+ \varepsilon h(x') + S_h(x',\varepsilon))-f(x',c)\Big)\sqrt{1+ \varepsilon^2|\nabla_{x'}h(x')|^2 + |\nabla_{x'}S_h(x', \varepsilon)|^2}\, dx'\\
&\quad +\lim_{\varepsilon\to 0}\int_{\mathbb{R}^{n-1} \cap Q^{n-1}_R(0)}f(x', c)\left(\frac{\sqrt{1+ \varepsilon^2|\nabla_{x'}h(x')|^2+ |c(\varepsilon)\nabla_{x'}s_h(x')|^2}-1}{\varepsilon}\right)dx'.    
\end{align*}
If we assume that our original function $h \in C^1_c(\mathbb{R}^{n-1})$ satisfies $\int_{\mathbb{R}^{n-1}}f(x', c)h(x')dx' = 0$, then $c(0) =0$ and \eqref{s_prime} imply the second limit vanishes.  Continuing under this assumption, we see that since $G_{\varepsilon} \Delta H(\vec{e_n}, c)$ is compact we may calculate
\begin{align*}
0 &=\lim_{\varepsilon\to 0}\frac{1}{\varepsilon}\int_{\mathbb{R}^{n-1} \cap Q^{n-1}_R(0)}\left(\int_c^{c+\varepsilon h(x')+S_h(x',\varepsilon)}\partial_{x_n}f(x',t)dt\right)\sqrt{1+ \varepsilon^2|\nabla_{x'}h(x')|^2+ |c(\varepsilon)\nabla_{x'}s_h(x')|^2}\, dx'\\
&=\lim_{\varepsilon\to 0}\int_{\mathbb{R}^{n-1} \cap Q^{n-1}_R(0)}\left(\frac{\varepsilon h(x')+S_h(x', \varepsilon)}{\varepsilon}\right)\left(\frac{\int_c^{c+\varepsilon h(x')+S_h(x', \varepsilon)}\partial_{x_n}f(x',t)dt}{\varepsilon h(x')+S_h(x', \varepsilon)}\right) \\
&\quad\cdot\sqrt{1+ \varepsilon^2|\nabla_{x'}h(x')|^2+ |c(\varepsilon)\nabla_{x'}s_h(x')|^2}\, dx'\\
&= \int_{\mathbb{R}^{n-1}}\partial_{x_n}f(x',c)h(x')dx'
\end{align*}
for $\mathcal{H}^1$-a.e. $c\in \mathbb{R}$ by the Lebesgue differentiation theorem (see \cite{evansgariepy}, Section 1.7, Corollary 2).

Therefore, for $\mathcal{H}^1$-a.e. $c\in \mathbb{R}$, for every $h \in C^1_c(\mathbb{R}^{n-1})$ such that $\int_{\mathbb{R}^{n-1}}f(x', c)h(x')dx' = 0$, then 
\begin{align*}
    \int_{\mathbb{R}^{n-1}}\partial_{x_n}f(x',c)h(x')dx' = 0.
\end{align*}
Thus, there must be a constant $K(c)$ such that for almost every $x' \in \mathbb{R}^{n-1}$,
\begin{align}\label{e:separation of variables}
\partial_{x_n}f(x',c) = K(c)f(x', c).
\end{align} 

\noindent {\bf Step 2}: Now, we claim that there exists a $\mathcal{H}^1$-measurable function $\overline{K}$ that satisfies 
\begin{align}\label{e:separation of variables 2}
\partial_{x_n}f(x',c) = \overline{K}(c)f(x', c)\mbox{\quad for $\mathcal{L}^n$-a.e. $(x',c)\in \mathbb{R}^n$}.
\end{align} 
To see the first claim, let
\begin{align*}
    G:=\{c \in \mathbb{R}: K(c) \text{ exists, \eqref{e:separation of variables} holds }\mathcal{H}^{n-1}\text{-a.e.}\}.
\end{align*}
We decompose $G$ as follows.
\begin{align*}
    G_p:= G \cap \{c \in \mathbb{R}: \mathcal{H}^{n-1}(\{x': f(x', c)\not =0\})>0\},\\
    G_z:= G \cap \{c \in \mathbb{R}: \mathcal{H}^{n-1}(\{x': f(x', c)\not =0\})=0\}.
\end{align*}
Note that $G$ is a set of full $\mathcal{H}^1$-measure, i.e., $\mathcal{H}^1(\mathbb{R}\setminus G)=0$, $G_p$ is $\mathcal{H}^1$-measurable by Fubini's theorem and the $\mathcal{H}^n$-measurability of $f$, and $G_z$ is $\mathcal{H}^1$-measurable by the $\mathcal{H}^1$-measurability of $G$ and $G_p$. We define the function $\widetilde{K}:\mathbb{R}^{n-1} \times \mathbb{R} \rightarrow \mathbb{R}$ as follows.  
\begin{align*}
    \widetilde{K}(x', c) := \begin{cases}
    \frac{\partial_{x_n}f(x'', c)}{f(x'', c)} & \mbox{if }(x',c) \in \mathbb{R}^{n-1}\times G_p, \\
    0& \mbox{if }(x',c) \in  \mathbb{R}^{n-1}\times G_p^c,
    \end{cases}
\end{align*}
where we interpret $\frac{\partial_{x_n}f(x'', c)}{f(x'', c)}$ to be the unique value that this ratio takes $\mathcal{H}^{n-1}$-almost everywhere on the set $\{(x', c): f(x', c)\not = 0\}$, for $c \in G_p$.

Note that $\widetilde{K}(x', c)$ is independent of $x' \in \mathbb{R}^{n-1}$.  We set $\overline{K}(c):= \widetilde{K}(x', c)$.

Now, we prove that $\overline{K}$ is $\mathcal{H}^1$-measurable.  Let $N:= \{f =0\}$.  Since $\partial_{x_n}f$ and $f$ are $\mathcal{H}^n$-measurable functions, $\widetilde{K}\big|_{ (\mathbb{R}^n \setminus N)}$ is $\mathcal{H}^n$-measurable. Recall that for any $k \in \mathbb{N}$, $\mathcal{H}^k$-measurability in $\mathbb{R}^k$ is equivalent to being approximately continuous $\mathcal{H}^k$-almost everywhere in $\mathbb{R}^k$ (see \cite{evansgariepy} Section 1.7.2, Theorem 3, and \cite{Federer69} Section 2.9.13). Therefore, we define for $c \in G_p$
\begin{align*}
    G_p(c):=\Big\{(x', c):\ & x' \in \mathbb{R}^{n-1}, (x', c) \text{ is a point of $\mathcal{L}^n$-density $1$ for the set $\mathbb{R}^n \setminus N$, and}\\
    & (x', c) \text{ is a point of approximate continuity of $\widetilde{K} \big|_{ (\mathbb{R}^n \setminus N)}$}\Big\}.
\end{align*}
Define $\widetilde{G}_p \subset G_p$ to be the set for which $G_p(c)$ is a set of full $\mathcal{H}^{n-1}$-measure in $\{(x', c): x' \in \mathbb{R}^{n-1}\} \cap (\mathbb{R}^n \setminus N)$. Note that this set $\widetilde{G}_p$ is a set of full $\mathcal{H}^1$-measure in $G_p$, and for $\mathcal{H}^{n-1}$-almost every $x' \in \{(x', c): x' \in \mathbb{R}^{n-1}\} \cap (\mathbb{R}^n \setminus N)$, $\widetilde{K}(x', c)= K(c)$.

For any $c_0 \in \widetilde{G}_p$ and any $(x'_0, c_0) \in G_p(c_0)$, we may calculate the following two equalities. Let $r>0$ and consider $Q_r^n(x'_0,c_0)$. First, by the assumption of density, for any $\varepsilon>0$,
\begin{align}\label{e:density assump}
    &\lim_{r \rightarrow 0} \frac{\mathcal{H}^n(Q_r^n(x'_0, c_0) \cap \{(x', c): |\widetilde{K}(x', c) - \widetilde{K}(x'_0, c_0)| \ge \varepsilon \})}{\mathcal{H}^n(Q_r^n(x'_0, c_0))} \\ \nonumber
    &  = \lim_{r \rightarrow 0} \frac{\mathcal{H}^n(Q_r^n(x'_0, c_0) \cap \{(x', c): |\widetilde{K}(x', c) - \widetilde{K}(x'_0, c_0)| \ge \varepsilon \} \setminus N )}{\mathcal{H}^n(Q_r^n(x'_0, c_0))}\\ \nonumber
    &  = \lim_{r \rightarrow 0} \frac{\int_{c_0-r/2}^{c_0 + r/2}\int_{Q_r^{n-1}(x'_0, c_0)} \chi_{\{|\widetilde{K}(x',c) - \widetilde{K}(x'_0, c_0)| \ge \varepsilon \}\setminus N}(x', c)d\mathcal{H}^{n-1}(x')d\mathcal{H}^1(c)}{\mathcal{H}^n(Q_r^n(x'_0, c_0))}\\ \nonumber
    &  = \lim_{r \rightarrow 0} \frac{\int_{c_0-r/2}^{c_0 + r/2} \chi_{\{|\overline{K}(c) - \overline{K}(c_0)| \ge \varepsilon \}}d\mathcal{H}^1(c)}{r} \\
	&\quad- \frac{\int_{c_0-r/2}^{c_0 + r/2}\int_{Q_r^{n-1}(x'_0, c_0)} \chi_{\{(x', c):|\overline{K}(c) - \overline{K}(c_0)| \ge \varepsilon \} \cap N}d\mathcal{H}^{n-1}(x')d\mathcal{H}^1(c)}{\mathcal{H}^n(Q_r^n(x'_0, c_0))}\nonumber\\ 
    &  = \lim_{r \rightarrow 0}\frac{\int_{c_0-r/2}^{c_0 + r/2} \chi_{\{|\overline{K}(c) - \overline{K}(c_0)| \ge \varepsilon \}}d\mathcal{H}^1(c)}{r}.\label{e: approx cont in 1D}
\end{align}

Second, by assumption of approximate continuity, for any $\varepsilon>0$, 
\begin{align}\label{e:approx cont assump}
\lim_{r \rightarrow 0} \frac{\mathcal{H}^n(Q^n_r(x'_0, c_0) \cap \{(x', c): |\widetilde{K}(x', c) - \widetilde{K}(x'_0, c_0)| \ge \varepsilon \} \setminus N )}{\mathcal{H}^n(Q^n_r(x'_0, c_0))}=0.
\end{align}

Since this holds for all $\varepsilon>0$, \eqref{e:density assump}, \eqref{e: approx cont in 1D}, and \eqref{e:approx cont assump} imply that $\overline{K}: \mathbb{R} \rightarrow \mathbb{R}$ is approximately continuous at all $c_0 \in \widetilde{G}_p$.  

It remains to show that $\overline{K}$ is approximately continuous $\mathcal{H}^1$-a.e. on $G_z$.  Define
\begin{align*}
    G_z(c):=\Big\{(x', c):\ & x' \in \mathbb{R}^{n-1}, (x', c) \text{ is a point of $\mathcal{L}^n$-density $1$ for the set $N$, and } \text{$\partial_{x_n}f(x', c)=0$}\Big\}.
\end{align*}
Define $\widetilde{G}_z \subset G_z$ to be the collection of $c \in G_z$ such that $G_z(c)$ is a set of full $\mathcal{H}^{n-1}$-measure. Note that $\mathcal{H}^{1}(G_z \setminus \widetilde{G}_z) =0$. Recall that for any $c \in \widetilde{G}_z$ we redefine $\overline{K}(c) = 0$ so that \eqref{e:separation of variables} holds. Now, by an identical chain of reasoning as in \eqref{e:density assump} and \eqref{e: approx cont in 1D}, we obtain
\begin{align*}
    0= & \lim_{r \rightarrow 0} \frac{\mathcal{H}^n(Q_r(x'_0, c_0) \cap \{(x', c): |\widetilde{K}(x', c) - \widetilde{K}(x'_0, c_0)| \ge \varepsilon \})}{\mathcal{H}^n(Q_r(x'_0, c_0))}\\ = & \lim_{r \rightarrow 0} \frac{\mathcal{H}^n(Q_r(x'_0, c_0) \cap \{(x', c): |\widetilde{K}(x', c) - \widetilde{K}(x'_0, c_0)| \ge \varepsilon \} \cap N )}{\mathcal{H}^n(Q_r(x'_0, c_0))}\\
     = & \lim_{r\rightarrow 0}\frac{\int_{c_0-r/2}^{c_0 + r/2} \chi_{\{|\overline{K}(c) - \overline{K}(c_0)| \ge \varepsilon \}}d\mathcal{H}^1(c)}{r}.
\end{align*}

Since $\overline{K}$ is approximately continuous at all points at $\widetilde{G}_p$ and $\widetilde{G}_z$, both $\widetilde{G}_p$ and $\widetilde{G}_z$ are $\mathcal{H}^1$-measurable, and $\mathcal{H}^1(\mathbb{R} \setminus (\widetilde{G}_p \cup \widetilde{G}_z))=0$, we conclude that $\overline{K}:\mathbb{R} \rightarrow \mathbb{R}$ is $\mathcal{H}^1$-measurable. Moreover,
$$\partial_{x_n}f(x',c) = \overline{K}(c)f(x', c)\mbox{\quad for $\mathcal{L}^n$-a.e. $(x',c)\in \mathbb{R}^n$}.$$
\noindent {\bf Step 3}: Now, we produce the collection of disjoint open intervals $\{I_j\}_j$ and show that $\overline{K}$ is locally integrable inside every $I_j$. Let $L_{x'} = \pi_{\mathbb{R}^{n-1}}^{-1}(x')$ and recall that for $\mathcal{H}^{n-1}$-a.e. $x'$ the function $f|_{ L_{x'}} \in W^{1, 1}_{\loc}(L_{x'}) \subset C(L_{x'})$. For such $x'$, $\{f |_{L_{x'}} >0\}$ is open in $\mathbb{R}$.  Thus, we may write 
\begin{align}\label{positive decomposition}
\left\{f |_{L_{x'}}>0\right\} = \bigcup^{\infty}_{j=1}\{x'\} \times I_j(x')
\end{align}
for some disjoint collection of open intervals $I_j(x') = (a_j(x'), b_j(x'))$. In particular, the functions $\partial_{x_n}f(x',\cdot)$ and $f(x', \cdot)$ are locally integrable in each $I_j(x')$, and hence $\overline{K}(\cdot)$ is locally integrable, as well. 

\noindent {\bf Step 4}: Now, we do the separation of variables.  For simplicity of notation, we shall relabel $\overline{K} =K$.  For any $x_{n, 0}, x_n \in I_j(x')$ we may integrate \eqref{e:separation of variables 2} to see that 
\begin{align}\label{FTC}
\log f(x',c)\big|_{c=x_{n, 0}}^{x_n}=\int_{x_{n, 0}}^{x_n}K(c)dc\implies f(x',x_n) =f(x',x_{n, 0})e^{\int_{x_{n,0}}^{x_n}K(c)dc}.
\end{align}
Let $Z$ be a set of $\mathcal{H}^{n-1}$-measure zero such that \eqref{positive decomposition} holds for all $x' \in \mathbb{R}^{n-1} \setminus Z$. Define
$$S=\left\{I_j(x'):x'\in \mathbb{R}^{n-1}\setminus Z,\ j\in \N\right\}.$$
We claim that $S$ is a collection of disjoint intervals of positive length and, therefore, $S$ is countable. Suppose not, there exists $I_j(x')\not= I_k(x'')$ such that $I_j(x')\cap I_k(x'')\not=\emptyset$. Without loss of generality, we may assume that $b_j(x')\in I_k(x'')= (a_k(x''), b_k(x''))$. Letting $x_{n} \rightarrow b_j(x')$ and $x_{n,0}\in I_j(x')\cap I_k(x'')$ in \eqref{FTC}, we observe that $\int_{x_{n,0}}^{b_j(x')} K(c)dc = -\infty$ since $f(x',b_j(x'))=0$ and $f(x',x_{n, 0})>0$. Moreover, $\int_{x_{n,0}}^{b_j(x')} K(c)dc = -\infty$ implies that $f(x'', b_j(x')) =0$ and hence $b_j(x')=b_k(x'')$. A similar argument shows that $a_j(x')=a_k(x'')$, leading to a contradiction. This means that we may find a collection of disjoint open intervals $S=\{I_j\}_j$ such that for $\mathcal{H}^{n-1}$-a.e. $x' \in \mathbb{R}^{n-1}$,
$$\{f(x', \cdot)>0 \} = \bigcup_{j \in A}I_j \quad\text{for some index set\ } A\subset \mathbb{N} \mbox{ depending on $x'$}.$$

Therefore, for each $I_j:=(a_j,b_j)$, we may choose $x_{n,0} = (a_j + b_j)/2$ in \eqref{FTC} and obtain that for every $I_j$, there exist functions $A_j(x'):=f(x',(a_j+b_j)/2)$ and $B_j(x_n):=e^{\int_{(a_j+b_j)/2}^{x_n}K(c)dc}$ such that 
\begin{align*}
    f(x', x_n) = A_j(x')B_j(x_n)
\end{align*}
for $\mathcal{H}^n$-a.e. $(x', x_n) \in \mathbb{R}^{n-1} \times I_j$.  This proves the lemma.
\end{proof}

\begin{remark}
    Lemma \ref{l:product structure} is a generalization of \cite{BrockChiacchio16} Theorem 2.1, in which the authors prove (among other things) that in a Lipschitz domain $\Omega\subset \mathbb{R}^n$, if $\nu = \phi\mathcal{L}^{n}$ for $\phi \in C^1(\Omega)$ and $\phi>0$ and the relative half spaces $\Omega \cap H(\vec{e_n}, c)$ are minimizers of $\operatorname{Per}_{\nu}(\cdot, \Omega)$ in their mass class for all $c \in \mathbb{R}$ then $\phi(x', x_n) = \rho(x')s(x_n)$ for $\rho, s$ both $C^1$ functions.  The assumption that $\phi \in C^1$ and $\phi>0$ allows the authors to use the Implicit Function theorem to obtain $s \in C^2$, and to force $\{I_j\}_j = \{ \mathbb{R}\}$.  However, under the more general assumptions of Lemma \ref{l:product structure} the possibility of multiple intervals $\{I_j\}$ is unavoidable. It is clear that Lemma \ref{l:product structure}, above, may be reformatted to handle such relative isoperimetric problems as well.
\end{remark}

Combining Lemma \ref{l:product structure} with Lemma \ref{l:all half-spaces} gives the following corollary.

\begin{corollary}\label{c:product structure in all directions}
Let $\mu \in \mathscr{W}(\mathbb{R}^n)$ with distribution function $f$.  Assume that $\mu$ a PS measure.  Then either $f \equiv 0$ or for every orthonormal frame $\{\vec{e_i} \}_i^n$ and coordinates $(x_1, ..., x_n)$ on $\mathbb{R}^n$ imposed by $\{\vec{e_i}\}_i^n$ \begin{align}\label{e:product formula}
        f\Big(\sum_{i=1}^nx_ie_i\Big) = \prod_{i=1}^nf_i(x_i)
    \end{align} for some functions $f_i \in L^1(\mathbb{R}) \cap W^{1, 1}_{\loc}(\mathbb{R})$ satisfying $f_i>0$, depending upon $\{\vec{e_i} \}_i^n$. In particular, we may choose the representative of $f$ defined by \eqref{e:product formula} everywhere and obtain that $f \in C(\mathbb{R}^n)$ and $f>0$. 
\end{corollary}

\begin{proof}
Without loss of generality, we may assume that $f\not \equiv 0$ and $\{\vec{e_k}\}_{k=1}^n$ is the standard basis of $\mathbb{R}^n$.  Lemma \ref{l:all half-spaces} allows us to apply Lemma \ref{l:product structure} in each direction $\vec{e_k}$ to produce a countable collection of axis-parallel rectangles $\{R_{j_1, ..., j_n}\}$ where 
\begin{align*}
    R_{j_1, ..., j_n} = I_{j_1} \times ... \times I_{j_n}
\end{align*}
for $I_{j_k} \in \{I^k_j\}_{j\in \N}$ the collection of intervals guaranteed by Lemma \ref{l:product structure} with $x_k$ in place of $x_n$.  

Now, suppose that we have an axis-parallel rectangle $R_{j_1, ..., j_n} $ in $\mathbb{R}^n$ and $n$ different decompositions of a function $f(x_1, ..., x_n) = A_i(x_1, ..., \hat{x_i}, ..., x_n)B_i(x_i)$ where we use the $\hat{x_i}$ notation to indicate that $A_i$ does not depend upon $x_i$. Then, by separation of variables, we may easily obtain that $f(x_1, ..., x_n) = \prod_{i=1}^n B_i(x_i)$. Thus, within each $R_{j_1, ..., j_n}$ there exist positive functions $f_{1}, ..., f_{n}$ with $f_k\in W^{1,1}_{}(I_{j_k})$ such that for $\mathcal{H}^{n}$-a.e. $(x_1, ..., x_n) \in R_{j_1, ..., j_n}$,
\begin{align*}
    f(x_1, ..., x_n) = \prod_{i=1}^n f_i(x_i).
\end{align*}
Thus, $f$ is positive and continuous on $R_{j_1, ..., j_n}$.

We claim that $\{I_{j}^k\}_{j\in \N} = \{\mathbb{R}\}$ for each $k =1, ..., n$ implying that there is only one rectangle, namely $\mathbb{R}^n$ itself.  To see the claim, observe that in general the zero set $\{f=0\}$ is of the form
\begin{align*}
    \{f=0\} = \bigcup_{i=1}^n\left(\mathbb{R}^{i-1} \times \{f_i=0\} \times \mathbb{R}^{n-i}\right).
\end{align*}
Note that this set is the (possibly uncountable) union of axis-parallel slabs of the form
\begin{align*}
   K= \mathbb{R}^{j-1} \times V \times \mathbb{R}^{n-j}
\end{align*}
for some connected, closed sets $V \subset \mathbb{R}$. If $\{I_{j}^k\}_{j\in \N}  \not = \{\mathbb{R}\}$ then $\{f=0\}$ is non-empty and there exists at least one set $K \subset \{f=0\} \not=\mathbb{R}^n$.  Hence, we may find a new orthonormal frame $(\varepsilon_1, ...,\varepsilon_{n})$ such that sets of the form $K$ are no longer axis-parallel.  But, because $f$ admits a product structure in this new basis, $\{f=0\}$ must have consist in a union of axis-parallel slabs in the new coordinate system.  If $f \not\equiv 0$, this provides a contradiction. Therefore, if $f\not \equiv 0$ we must have $\{f=0 \} = \emptyset$.  In particular, $R_{j_1, ..., j_n} = \mathbb{R}^n$ and the claim holds.

Thus, we may choose the function defined by \eqref{e:product formula} to be the representative of $f \in L^1(\mathbb{R}^n) \cap W^{1, 1}_{\loc}(\mathbb{R}^n)$. The rest of the properties follow immediately from this choice and the fact that the $f_i \in W^{1,1}_{\loc}(\mathbb{R}) \subset C(\mathbb{R})$.
\end{proof}

\section{Proof of Theorem \ref{t:main theorem 2}}\label{s:symm and proof of main theorem}

\noindent {\bf Step 1}: Since $\mu\in \mathscr{W}(\mathbb{R}^n)$ is a PS measure, by Corollary \ref{c:product structure in all directions}, if $f\equiv 0$, we can choose $C=0$. If $f\not\equiv 0$, our objective is to demonstrate the existence of constants $0< c< \infty$, $0 < C < \infty$, and $a\in \mathbb{R}^n$ such that
\begin{align*}
    f(x) = Ce^{-c|x-a|^2}.
\end{align*}
With the aid of Corollary \ref{c:product structure in all directions}, we can write
$$f(x_1,...,x_n)=f_1(x_1)\cdots f_n(x_n).$$
Applying Lemma \ref{l:Ehrhard implies symmetric_3} to $f$, there exists $m\in \mathbb{R}$ such that 
$$\nu=\widetilde{f}\mathcal{L}^n$$
is a PS measure, where $\widetilde{f}(x)=f(x+m)$, $\widetilde{f}(x_1,...,x_n)=\widetilde{f}_1(x_1)\cdots \widetilde{f}_n(x_n)$, $\widetilde{f}_i>0$, and $\widetilde{f}_i(-t)=\widetilde{f}_i(t)$ for all $i=1,...,n$. Let $a=m\in \mathbb{R}^n$. Our aim is to establish that
$$\widetilde{f}(x)=Ce^{-c|x|^2}$$
for some $0< c< \infty$, $0 < C < \infty$. For ease of reference, we will continue to use $f$ to denote $\widetilde{f}$ throughout the remainder of the proof.\\
\noindent {\bf Step 2}: In Step 1, we established that $f(x_1,...,x_n)=f_1(x_1)\cdots f_n(x_n)$, $f_i\in C(\mathbb{R})$, $f_i>0$ and $f_i(-t)=f_i(t)$ for all $i=1,...,n$. Consequently, each component function $f_i$ can be expressed as
\begin{align}\label{exp_sum}
f_i(t)=f_i(0)\left(\frac{f_i(t)}{f_i(0)}\right)=c_ie^{g_i(t)}
\end{align}
for some continuous function $g_i$ and $g_i(-t)=g_i(t)$, where $c_i=f_i(0)$ and $g_i(0)=0$. That is,
$$f(x_1,...,x_n)=Ce^{g_1(x_1)+\cdots+g_n(x_n)},$$
where $C:=f(0,...,0)=f_1(0)\cdots f_n(0)>0$. We claim that
\begin{align}\label{g_1=g_2}
g_1=g_2=\cdots =g_n.
\end{align} 
Consider
$$O=\left[ \begin{array}{cc} O_2 & \mathbf{0} \\ \mathbf{0} & I_{n-2} \end{array} \right]\in M_n(\mathbb{R}),\quad O_2=\left[ \begin{array}{cc} \frac{1}{\sqrt{2}} & \frac{1}{\sqrt{2}} \\ \frac{1}{\sqrt{2}}& \frac{-1}{\sqrt{2}}\end{array} \right]\in M_2(\mathbb{R}).$$
Then
\begin{align}\label{exp_with_g}
f(Ox)&=f\left(\frac{1}{\sqrt{2}}x_1+\frac{1}{\sqrt{2}}x_2, \frac{1}{\sqrt{2}}x_1-\frac{1}{\sqrt{2}}x_2,x_3,\cdots x_n\right)\notag\\
&=C\exp\left(g_1\left(\frac{1}{\sqrt{2}}x_1+\frac{1}{\sqrt{2}}x_2\right)+g_2\left(\frac{1}{\sqrt{2}}x_1-\frac{1}{\sqrt{2}}x_2\right)+\cdots+ g_n(x_n)\right).
\end{align}
Utilizing Lemma \ref{c:product structure in all directions} with the orthogonal frame $\{Oe_i\}_{i=1}^n$ and the argument in \eqref{exp_sum}, we obtain
\begin{align}\label{exp_with_h}
f(Ox)=f\left(\sum_{i=1}^nx_iOe_i\right)=Ce^{h_1(x_1)+\cdots+h_n(x_n)},
\end{align}
where $h_i$ depends on matrix $O$ and satisfies $h_i(0) = 0$. To confirm the property of \eqref{g_1=g_2}, we first check that
\begin{align}\label{property_of_h}
h_1(\alpha)=h_1(-\alpha),\quad h_1(\alpha)=h_2(\alpha),
\end{align}
for all $\alpha\in \mathbb{R}$. Setting $x_3=x_4=\cdots =x_n=0$ in \eqref{exp_with_g} and \eqref{exp_with_h}, we have
\begin{align}\label{h_1_plus_h_2}
h_1(x_1)+h_2(x_2)=g_1\left(\frac{1}{\sqrt{2}}x_1+\frac{1}{\sqrt{2}}x_2\right)+g_2\left(\frac{1}{\sqrt{2}}x_1-\frac{1}{\sqrt{2}}x_2\right).
\end{align}
In particular, letting $x_1 = \alpha$ and $x_2 = 0$, we find
\begin{align*}
h_1(\alpha)=g_1\left(\frac{1}{\sqrt{2}}\alpha\right)+g_2\left(\frac{1}{\sqrt{2}}\alpha\right)=h_1(-\alpha),
\end{align*}
where we have used the fact that $g_i(-t)=g_i(t)$. Similarly, for $x_1 = 0$ and $x_2 = \alpha$, we have
\begin{align*}
h_2(\alpha)=g_1\left(\frac{1}{\sqrt{2}}\alpha\right)+g_2\left(-\frac{1}{\sqrt{2}}\alpha\right)=g_1\left(\frac{1}{\sqrt{2}}\alpha\right)+g_2\left(\frac{1}{\sqrt{2}}\alpha\right)=h_1(\alpha).
\end{align*}
Consider $x_1=x_2=\frac{\alpha}{\sqrt{2}}$ in \eqref{h_1_plus_h_2} and apply \eqref{property_of_h}, we also have
$$g_1(\alpha)=h_1\left(\frac{\alpha}{\sqrt{2}}\right)+h_2\left(\frac{\alpha}{\sqrt{2}}\right)=2h_1\left(\frac{\alpha}{\sqrt{2}}\right).$$
Similarly, for $x_1 = \frac{\alpha}{\sqrt{2}}$ and $x_2 = -\frac{\alpha}{\sqrt{2}}$, we apply \eqref{property_of_h} to obtain,
$$g_2(\alpha)=h_1\left(\frac{\alpha}{\sqrt{2}}\right)+h_2\left(-\frac{\alpha}{\sqrt{2}}\right)=h_1\left(\frac{\alpha}{\sqrt{2}}\right)+h_1\left(-\frac{\alpha}{\sqrt{2}}\right)=2h_1\left(\frac{\alpha}{\sqrt{2}}\right).$$
Hence, $g_1 = g_2$. An identical argument yields $g_1 = g_2 = \cdots = g_n$.\\
{\bf Step 3}: We claim that for all positive integers $k$, we have
\begin{align}\label{recursive equation}
    2g_1(\sqrt{k}\alpha)+2g_1(\alpha)=g_1((\sqrt{k}-1)\alpha)+g_1((\sqrt{k}+1)\alpha).
\end{align}
Consider
$$O=\left[ \begin{array}{cc} O_2 & \mathbf{0} \\ \mathbf{0} & I_{n-2} \end{array} \right]\in M_n(\mathbb{R}),\quad O_2=\left[ \begin{array}{cc} \frac{\sqrt{k+k^2}}{k+1} & -\frac{\sqrt{k+1}}{k+1} \\  \frac{\sqrt{k+1}}{k+1}&\frac{\sqrt{k+k^2}}{k+1} \end{array} \right]\in M_2(\mathbb{R}).$$
Using the argument in Step 2 and \eqref{g_1=g_2}, we have
\begin{align*}
f(Ox)&=C\exp\left(g_1\left(\frac{\sqrt{k+k^2}}{k+1}x_1-\frac{\sqrt{k+1}}{k+1}x_2\right)+g_1\left(\frac{\sqrt{k+k^2}}{k+1}x_1+\frac{\sqrt{k+1}}{k+1}x_2\right)+\cdots+ g_1(x_n)\right),
\end{align*}
and
\begin{align*}
f(Ox)=f\left(\sum_{i=1}^nx_iOe_i\right)=Ce^{\ell_1(x_1)+\cdots+\ell_n(x_n)},
\end{align*}
for some $\ell_i$ depending on matrix $O$ and $\ell_i(0)=0$. Setting $x_3=x_4=\cdots =x_n=0$,
\begin{align*}
\ell_1(x_1)+\ell_2(x_2)=g_1\left(\frac{\sqrt{k+k^2}}{k+1}x_1-\frac{\sqrt{k+1}}{k+1}x_2\right)+g_1\left(\frac{\sqrt{k+k^2}}{k+1}x_1+\frac{\sqrt{k+1}}{k+1}x_2\right).
\end{align*}
In particular, setting $x_2 = 0$, we find
\begin{align*}
\ell_1(x_1)=2g_1\left(\frac{\sqrt{k}}{\sqrt{k+1}}x_1\right).
\end{align*}
Similarly, when $x_1 = 0$, we get
\begin{align*}
\ell_2(x_2)=2g_1\left(\frac{1}{\sqrt{k+1}}x_2\right),
\end{align*}
where we have again utilized the fact that $g_i(-t)=g_i(t)$. Thus,
\begin{align*}
&2g_1\left(\frac{\sqrt{k}}{\sqrt{k+1}}x_1\right)+2g_1\left(\frac{1}{\sqrt{k+1}}x_2\right)\\
&=g_1\left(\frac{\sqrt{k+k^2}}{k+1}x_1-\frac{\sqrt{k+1}}{k+1}x_2\right)+g_1\left(\frac{\sqrt{k+k^2}}{k+1}x_1+\frac{\sqrt{k+1}}{k+1}x_2\right).
\end{align*}
Finally, if we set $x_1 = x_2 = \sqrt{k+1}\alpha$, we conclude that
$$2g_1(\sqrt{k}\alpha)+2g_1(\alpha)=g_1((\sqrt{k}-1)\alpha)+g_1((\sqrt{k}+1)\alpha).$$
{\bf Step 4}: We claim that for all positive integers $k$, we have
\begin{align}\label{e: base relation sqrt}
    g_1(k\alpha)=k^2g_1(\alpha).
\end{align}
We will prove this by induction. Note that the base case $k=1$ is trivial. Now, assume that the claim holds for all integers less than or equal to $k$. For any $m\in \mathbb{N}$, substituting $k=m^2$ into \eqref{recursive equation}, we obtain
\begin{align*}
2g_1(m\alpha)+2g_1(\alpha)=g_1((m-1)\alpha)+g_1((m+1)\alpha).
\end{align*}
That is,
\begin{align}\label{induction_on_g}
g_1(m\alpha)-g_1((m-1)\alpha)+2g_1(\alpha)=g_1((m+1)\alpha)-g_1(m\alpha).
\end{align}
Summing $m$ from $1$ to $k$ in \eqref{induction_on_g}, we have
\begin{align*}
g_1(k\alpha)+2kg_1(\alpha)&=\sum_{m=1}^k\Big(g_1(m\alpha)-g_1((m-1)\alpha)+2g_1(\alpha)\Big)\\
&=\sum_{k=1}^k\Big(g_1((m+1)\alpha)-g_1(m\alpha)\Big)=g_1((k+1)\alpha)-g_1(\alpha).
\end{align*}
Therefore,
$$g_1((k+1)\alpha)=g_1(k\alpha)+2kg_1(\alpha)+g_1(\alpha)=k^2g_1(\alpha)+2kg_1(\alpha)+g_1(\alpha)=(k+1)^2g_1(\alpha),$$
where we have used the induction hypothesis on $k$. This proves the claim. \\
\noindent {\bf Step 5}: Finally, we define
$$h(\alpha):=\frac{g_1(\alpha)}{\alpha^2}$$
for all $\alpha>0$. By \eqref{e: base relation sqrt}, we know $g_1(k\alpha)=k^2g_1(\alpha)$ and hence $h(k\alpha)=h(\alpha)$ for any positive integer $k$. In particular,
$$h(1)=h(1/q)$$
for all positive integers $q$. For any positive rational number $p/q$, we have
$$h(p/q)=h\left(p\left(\frac{1}{q}\right)\right)=h(1/q)=h(1).$$
For any $\alpha>0$, there exists a sequence of rational numbers $\{r_k\}$ such that $r_k\to \alpha$. Since $h$ is continuous at $\alpha$, we see that
$$h(1)=\lim_{k\to \infty}h(r_k)=h(\alpha)$$ for all $\alpha>0$. Letting $c=-h(1)$ and unwinding the definition of $h$, we have $g_1(\alpha)=-c\alpha^2$ for all $\alpha>0$. By $g_1(\alpha)=g_1(-\alpha)$ and $g_1(0)=0$,
$$g_1(\alpha)=-c\alpha^2$$
for all $a\in \mathbb{R}$. To complete the proof, we only need to prove that the constant $c$ is positive. However, if $c\le 0$, then $f \not \in L^1(\mathbb{R}^n)$.  Therefore, 
$$f(x)=Ce^{g_1(x_1)+\cdots+g_n(x_n)}=Ce^{-c(x_1^2+x_2^2+\cdots +x_n^2)}=Ce^{-c|x|^2}.$$
That is, $f$ is an isotropic Gaussian. \qed \\

\noindent {\bf Remark}: The constant $c>0$ also follows from \cite{Chambers19}, where Chambers proves that balls are the unique perimeter minimizers of measures with distributions that are smooth log-convex, rotationally symmetric functions. Given that Lemma \ref{l:H are minimal} establishes all half-spaces as $\mu$-perimeter minimizers, it follows that $c>0$. We also note that after Step 2 in the above proof, we have that $\mu = \mu_1 \times ... \times \mu_1$ is the $n$-fold product of some measure $\mu_1$ and half-spaces are the minimizers for perimeter.  In the context of the $\varepsilon$-enlargement definition of isoperimetric problems, Bobkov and Udre have shown that these two conditions imply that $\mu$ is an isotropic Gaussian \cite{BobkovUdre96_halfspaces_minimal_for_productmeasures_implies_Gaussian}.

\section{Appendix A: An example in one dimension}\label{s:appendix A examples}

As mentioned in the Introduction, Theorem \ref{t:main theorem 2} only works in dimensions $n \ge 2$. Inspired by the examples from \cite{Bobkov96} at the end of Section 1, we can similarly construct an example of a PS measure in $1$-dimension where the distribution does not have the form of a Gaussian distribution. The main ingredients for the proof are \cite{Bobkov1997} Proposition 13.4 and the argument in \cite{Yeh2023} Proposition 4.6 (3).

\begin{example}[An Example in $\mathbb{R}$]\label{ex:logconcave}
Let $\mu = f \mathcal{L}^1$, where the standard logistic distribution function $f$ is defined as
\begin{align*}
    f(x) = \displaystyle\frac{1}{(e^{\frac{x}{2}}+e^{-\frac{x}{2}})^2}.
\end{align*}
Then $\mu$ is a PS measure. However, $f$ does not have the form of a Gaussian distribution.
\end{example}
\begin{proof}
Let $E$ be a set of finite $\mu$-perimeter. Since $f$ is smooth and locally bounded away from zero, $E$ is also a set of locally finite $\mathcal{L}^1$-perimeter. Our goal is to show that
$$\operatorname{Per}_\mu(E^s_{-\vec{e}_1})\leq \operatorname{Per}_\mu(E)$$
where $\vec{e}_1=1$. Then, by the symmetry of $f$, $\operatorname{Per}_\mu(E^s_{\vec{e}_1})\leq \operatorname{Per}_\mu(E)$ and hence $\mu$ is a PS measure. Let
$$F(t)=\int_{-\infty}^t f(x)dx,\quad y=F^{-1}\left(\mu(E)\right).$$
Then $\mu(E)=\mu((-\infty,y])$ and hence
$$E^s_{-\vec{e_1}}=(-\infty,y].$$
{\bf Step 1}: Assume that 
$$E= \bigcup_{i=1}^m (a_i, b_i)$$ 
is the union of $m$ disjoint intervals such that $\{(a_i, b_i)+B_h\}_{i=1}^m$ are pairwise disjoint. We claim that
$$\mu(E+B_h)\geq \mu(E^s_{-\vec{e_1}}+B_h)\quad \text{for all  } h>0.$$
Since $F(x)=\mu((-\infty,x])$, we have $F(x)=\frac{1}{1+e^{-x}}$, $F^{-1}(p)=\log(\frac{p}{1-p})$, and
$$f(F^{-1}(p))=p(1-p).$$
Moreover,
\begin{align*}
f(F^{-1}(p+q))=(p+q)(1-p-q)&\leq p(1-p)-2pq+q(1-q)\\
&\leq f(F^{-1}(p))+f(F^{-1}(q)).
\end{align*}
Now we let $p=\mu(E)\in (0,1)$. Applying \cite{Bobkov1997} Proposition 13.4, 
$$\mu(E+B_h)\geq \min_{\mu(E)=p}\mu(E+B_h)=\mu(E^s_{-\vec{e_1}}+B_h)\mbox{\quad for all $h>0$.}$$

\noindent{\bf Step 2}: Since $E$ is a set of locally finite perimeter, $E\cap (-R,R)$ is a set of finite perimeter for any $R>0$. Applying \cite{Maggi_book12} Proposition 12.13 on $E\cap (-R,R)$, we may assume that $E\cap (-R,R)$ is a disjoint union of open intervals, pairwise separated by some positive distance. Since $f$ is locally bounded away from zero, $E\cap (-R,R)$ is a finite disjoint union of open intervals with positive distance. That is,
$$E\cap (-R,R)=\bigcup_{i\in S_R} (a_i, b_i)$$
where $|S_R|<\infty$. Applying Step 1 with $E\cap (-R,R)$, we have
$$\mu\left(E\cap (-R,R)+B_h\right)\geq \mu\left(\left(E\cap (-R,R)\right)_{-\vec{e_1}}^s+B_h\right).$$
By the definition of $\mu$ and the fundamental theorem of calculus,
\begin{align*}
\operatorname{Per}_\mu(E\cap (-R,R))&=\sum_{i\in S_R}\frac{1}{(e^{\frac{a_i}{2}}+e^{-\frac{a_i}{2}})^2}+\frac{1}{(e^{\frac{b_i}{2}}+e^{-\frac{b_i}{2}})^2}\\
&=\lim_{h\to 0^+}\sum_{i\in S_R}\frac{\mu((a_i,b_i)+B_h)-\mu(a_i,b_i)}{h}\\
&\geq\lim_{h\to 0^+}\frac{1}{h}\left[\mu\left(\bigcup_{i\in S_R}(a_i,b_i)+B_h\right)-\mu\left(\bigcup_{i\in S_R}(a_i,b_i)\right)\right]\\
&= \lim_{h\to 0^+}\frac{1}{h}\Big[\mu\left(E\cap (-R,R)+B_h\right)-\mu\left(E\cap (-R,R)\right)\Big]\\
&\geq  \lim_{h\to 0^+}\frac{1}{h}\Big[\mu\left(\left(E\cap (-R,R)\right)_{-\vec{e_1}}^s+B_h\right)-\mu\left(\left(E\cap (-R,R)\right)_{-\vec{e_1}}^s\right)\Big]\\
&=\frac{1}{(e^{\frac{y_R}{2}}+e^{-\frac{y_R}{2}})^2}
\end{align*}
where
$$y_R=F^{-1}\left(\mu(E\cap (-R,R))\right),\quad (E\cap (-R,R))_{\vec{e_1}}^s=(-\infty,y_R].$$
Taking $R\to \infty$, we have $y_R\to y$ and hence
$$\operatorname{Per}_\mu(E)\geq \frac{1}{(e^{\frac{y}{2}}+e^{-\frac{y}{2}})^2}=\operatorname{Per}_\mu((-\infty,y])=\operatorname{Per}_\mu(E^s_{-\vec{e_1}}).$$
\end{proof}

\section{Appendix B: Measurability}\label{s:app B measurability}

In this section we prove Theorem \ref{l:measurability}, which we restate below. 

\begin{theorem}[Theorem \ref{l:measurability}]
    Let $\mu \in \mathscr{W}(\mathbb{R}^n)$ with distribution function $f$. The generalized $\vec{v}$-Ehrhard symmetrization with respect to $\mu$ (see Definition \ref{d:gen Ehrh symm}) preserves measurability.  That is, for any $\mathcal{L}^n$-measurable set $E \subset \mathbb{R}^n$ and any direction $\vec{v} \in \mathbb{S}^{n-1}$, $S_{\vec{v}}(E)$ is a $\mathcal{L}^n$-measurable set.
\end{theorem}

By rotation, we may assume that $\vec{v} = -\vec{e_n}$.  Furthermore, we write $x= (x',x_n) \in \mathbb{R}^{n-1} \times \mathbb{R}$.  For any $E\subset \mathbb{R}^n$ and $x' \in \mathbb{R}^{n-1}$, we define the {\bf $x'$-slice} of $E$ by 
\begin{align*}
    E_{x'} = E \cap \{(x', x_n): x_n \in \mathbb{R}\}.
\end{align*}
For a function $f \in L^1(\mathbb{R}^n; [0,\infty])$, we define the {\bf marginal slice function of $E$} with respect to $f$ as 
\begin{align*}
    m_E(f, x') =  \int_{E_{x'}}f(x', x_n)dx_n.
\end{align*}
Note that $m_E(f, \cdot)$ is defined $\mathcal{L}^{n-1}$-almost everywhere. Let $\Omega$ be an open set in $\mathbb{R}^{n-1}$ and $v:\Omega \rightarrow [0,\infty]$ be a non-negative function.  A {\bf largest possible value function $h:\Omega\to [-\infty,\infty]$ of $v$ (with respect to $f$)} is a function that satisfies
\begin{enumerate}
\item[(1)] For $\mathcal{L}^{n-1}$-a.e. $x'$ in $\Omega$,
$$v(x') = \int_{-\infty}^{ h(x')}f(x', x_n)dx_n.$$
\item[(2)] For any function $\tilde{h}:\Omega\to [-\infty,\infty]$ with 
\begin{align*}
    v(x') & = \int_{-\infty}^{ \tilde{h}(x')}f(x', x_n)dx_n \mbox{\quad for $\mathcal{L}^{n-1}$-a.e. $x'$ in $\Omega$},
\end{align*}
we have $\tilde{h}(x')\leq h(x')$ for $\mathcal{L}^{n-1}$-a.e. $x'$ in $\Omega$.
\end{enumerate}
\begin{remark}
Let $S_{\vec{v}}(E)$ be the generalized $\vec{v}$-Ehrhard symmetrization with respect to $\mu=f\mathcal{L}^n$ with $\vec{v}=-\vec{e_n}$. Then
$$S_{\vec{v}}(E)=\{(x', x_n)\in \mathbb{R}^n: x_n \le h(x')\}$$ 
for some largest possible value function $h:\mathbb{R}^{n-1}\to [-\infty,\infty]$ of $m_E(f,\cdot)$.
\end{remark}

The proof of Theorem \ref{l:measurability} depends upon the following technical lemmata. A similar result for the case where the distribution function $f$ is an anisotropic Gaussian can be found in \cite{Yeh2023} Theorem 4.5.

\begin{lemma}\label{l:meas implies meas arbitrary I}
    Let $E$ be a $\mathcal{L}^n$-measurable set in $\mathbb{R}^n$ and $f\in L^1(\mathbb{R}^n)\cap W^{1,1}_{\loc}(\mathbb{R}^n)$ with $f\geq 0$. For any bounded open set $\Omega\subset \mathbb{R}^{n-1}$, there exists a $\mathcal{L}^{n-1}\res \Omega\, $-measurable function $h$ such that 
\begin{align*}
    m_E(f,x') & = \int_{-\infty}^{ h(x')}f(x', x_n)dx_n
\end{align*}
for $\mathcal{L}^{n-1}$-a.e. $x'$ in $\Omega$.
\end{lemma}

\begin{lemma}\label{l:meas implies meas arbitrary II}
        Let $E$ be a $\mathcal{L}^n$-measurable set in $\mathbb{R}^n$ and $f\in L^1(\mathbb{R}^n)\cap W^{1,1}_{\loc}(\mathbb{R}^n)$ with $f\geq 0$. Then there exists a $\mathcal{L}^{n-1}$-measurable largest possible value function $h:\mathbb{R}^{n-1}\to [-\infty,\infty]$ of $m_E(f,\cdot)$ (thought of as a ``heigh function'').
\end{lemma}

With Lemma \ref{l:meas implies meas arbitrary II}, we are able to prove Theorem \ref{l:measurability} as follows.

\subsection{Proof of Theorem \ref{l:measurability} (assuming Lemma \ref{l:meas implies meas arbitrary II})} 
Let $E$ be a measurable set in $\mathbb{R}^n$. Suppose 
$$S_{\vec{v}}(E)=\{(x', x_n)\in \mathbb{R}^n: x_n \le h(x')\}$$ 
for some largest possible value function $h:\mathbb{R}^{n-1}\to [-\infty,\infty]$ of $m_E(f,\cdot)$. Our goal is to show that $h$ is $\mathcal{L}^{n-1}$-measurable. Applying Lemma \ref{l:meas implies meas arbitrary II} to $m_E(f,\cdot )$, there exists a $\mathcal{L}^{n-1}$-measurable largest possible value function $\tilde{h}:\mathbb{R}^{n-1}\to [-\infty,\infty]$ of $m_E(f,\cdot )$. Since $\tilde{h}$ and $h$ are both largest possible value functions of $m_E(f,\cdot)$, we have $h(x')=\tilde{h}(x')$ for $\mathcal{L}^{n-1}$-a.e. $x'$ in $\mathbb{R}^{n-1}$. That is, $h$ is $\mathcal{L}^{n-1}$-measurable and hence $S_{\vec{v}}(E)$ is $\mathcal{L}^n$-measurable (see, for example, the argument in \cite{evansgariepy}, Section 2.2, Lemma 1).
\qed

\subsection{Proof of Lemma \ref{l:meas implies meas arbitrary I}}
To prove Lemma \ref{l:meas implies meas arbitrary I}, we construct a sequence of ``height functions'' with respect to distribution functions $f_i$ which approximate $f$.\\ 

\noindent {\bf Step 1}: We show that there is a sequence of functions $f_i \in C^{\infty}(\mathbb{R}^n)\cap L^1(\mathbb{R}^n)$ such that 
\begin{enumerate}
    \item[(i)] $f_i \to f \quad \text{in } L^1(\mathbb{R}^n)$.
    \item[(ii)] $f_i(x) >0 \quad \forall x \in \mathbb{R}^{n}$.
    \item[(iii)] For each $i \in \mathbb{N}$, the marginal function 
    \begin{align*}
        m_{\mathbb{R}^n}(f_i, x'):=\int_{\mathbb{R}} f_i(x', x_n) dx_n
    \end{align*}
    is continuous as a function of $x' \in \mathbb{R}^{n-1}$.  In fact, for any fixed $z \in \mathbb{R} \cup \{-\infty, \infty\}$ the function
    \begin{align*}
        \int_{-\infty}^{z} f_i(x', x_n) dx_n
    \end{align*}
    is continuous as a function of $x' \in \mathbb{R}^{n-1}$.
\end{enumerate}
Since $f\geq 0$, there exists $\tilde{f}_i\geq 0$ in $C_c^{\infty}(\mathbb{R}^{n})$ such that $\tilde{f}_i\to f$ in $L^1(\mathbb{R}^n)$. Define
$$f_i(x):=\tilde{f}_i(x)+\frac{1}{i}\frac{1}{(2\pi)^n}e^{-|x|^2/2}.$$
It is easy to see that $\{f_i\}_i$ satisfies (i) and (ii) since
$$\int_{\mathbb{R}^n}|f-f_i|dx\leq \int_{\mathbb{R}^n}|f-\tilde{f}_i|dx+\frac{1}{i}\int_{\mathbb{R}^n}\frac{1}{(2\pi)^n}e^{-|x|^2/2}dx\to 0.$$
Now, fix $i$ and $z$.  For any $x'_0\in \mathbb{R}^{n-1}$, we claim that
$$\lim_{x'\to x'_0}\int_{-\infty}^{z} f_i(x', x_n) dx_n=\int_{-\infty}^{z} f_i(x'_0, x_n) dx_n.$$
Note that there exists $R>0$ such that $x_0'\in (-R,R)^{n-1}$ and $\mbox{spt} |\nabla \tilde{f}_i|\subset [-R,R]^n$. Thus, 
\begin{align*}
&\left|\int_{-\infty}^{z} f_i(x', x_n) dx_n-\int_{-\infty}^{z} f_i(x_0', x_n) dx_n\right|\\
&\leq \int_{-\infty}^{z} \left|f_i(x', x_n) - f_i(x_0', x_n)\right| dx_n\\
&\leq \int_{-\infty}^{z} \left|\tilde{f}_i(x', x_n)- \tilde{f}_i(x_0', x_n)\right|dx_n+\frac1i\frac{1}{(2\pi)^n}\int_{-\infty}^{z}\left|e^{-|(x',x_n)|^2/2}-e^{-|(x'_0,x_n)|^2/2}\right| dx_n\\
&\leq \int_{-\infty}^{z} \left|\nabla \tilde{f}_i(\zeta(x_n),x_n) \right||x'-x'_0|\,dx_n+\frac1i\frac{1}{(2\pi)^n}\int_{-\infty}^{z}\left|e^{-|x'|^2/2}-e^{-|x'_0|^2/2}\right|e^{-x_n^2/2}\, dx_n\\
&\leq |x'-x_0'|(2R)\Big(\sup_{ [-R,R]^n}|\nabla \tilde{f}_i |\Big)+\frac1i\frac{1}{(2\pi)^n}\left|e^{-|x'|^2/2}-e^{-|x'_0|^2/2}\right|\int_{\mathbb{R}}e^{-x_n^2/2}\, dx_n\to 0
\end{align*}
as $x'\to x'_0$, where $\zeta(x_n)$ lies between $x'$ and $x'_0$. This establishes (iii).\\

\noindent {\bf Step 2}: For any $\mathcal{L}^n$-measurable set $E \subset \mathbb{R}^n$ and any $x'\in \Omega$, by Step 1 (ii), there exists a unique value $h_i(x')$ such that 
\begin{align*}
   m_E(f_i, x')= \int_{-\infty}^{h_i(x')}f_i(x', x_n)dx_n.
\end{align*}
Since $m_E(f_i, \cdot)$ is measurable and $\mathcal{L}^{n-1}(\Omega)<\infty$, by Lusin's theorem, for any $\varepsilon>0$ we may find a closed set $C_\varepsilon\subset \Omega$ and a continuous function $g^\varepsilon_i: \Omega \rightarrow [0,\infty)$ such that $\mathcal{L}^{n-1}(\Omega\setminus C_\varepsilon) < \varepsilon$ and $g^\varepsilon_i=m_E(f_i,\cdot)$ on $C_\varepsilon$. We claim that $h_i$ is continuous on $C_\varepsilon$. Let $x' \in C_\varepsilon$.  Suppose for the sake of contradiction that $\{ x_j\}_j \subset C_{\varepsilon}$ is a sequence such that $x_j \rightarrow x'$, $\lim\limits_{j\to \infty} h_i(x_j)$ exists (allowing for $\pm \infty$), and $h_i(x_j) \not\rightarrow h_i(x')$. Without loss of generality, we may assume that $h_i(x') > \lim\limits_{j\to \infty} h_i(x_j)$. Now we calculate the following:
\begin{align*}
g^\varepsilon_{i}(x')-g^\varepsilon_{i}(x'_j)&=\int_{-\infty}^{h_i(x')}f_i(x', x_n)dx_n-\int_{-\infty}^{h_i(x'_j)}f_i(x_j', x_n)dx_n\\
&=\int_{-\infty}^{h_i(x')}f_i(x', x_n)dx_n-\int_{-\infty}^{h_i(x')}f_i(x'_j, x_n)dx_n\\
&\quad+\int_{-\infty}^{h_i(x')}f_i(x'_j, x_n)dx_n-\int_{-\infty}^{h_i(x'_j)}f_i(x_j', x_n)dx_n\\
&=\int_{h_i(x_j')}^{h_i(x')}f_i(x_j', x_n)dx_n+\int_{-\infty}^{h_i(x')}f_i(x', x_n)-f_i(x'_j, x_n)dx_n.
\end{align*}
Since $h_i(x')$ is fixed, by Step 1 (iii), as $j\to \infty$,
$$\int_{-\infty}^{h_i(x')}f_i(x', x_n)-f_i(x'_j, x_n)dx_n\to 0.$$
By the continuity of $g_i^\varepsilon$, we have
\begin{align*}
0=\lim\limits_{j\to \infty}\int_{h_i(x_j')}^{h_i(x')}f_i(x_j', x_n)dx_n&\geq \lim\limits_{j\to \infty}\int_{h_i(x_j')}^{h_i(x')}\frac{1}{i}\frac{1}{(2\pi)^n}e^{-|(x_j',x_n)|^2/2}dx_n\\
&=\frac{1}{i}\frac{1}{(2\pi)^n}\lim\limits_{j\to \infty}e^{-|x_j'|^2/2}\int_{h_i(x_j')}^{h_i(x')}e^{-x_n^2/2}dx_n\\
&=\frac{1}{i}\frac{1}{(2\pi)^n}e^{-|x'|^2/2}\int_{\lim\limits_{j\to \infty}h_i(x_j')}^{h_i(x')}e^{-x_n^2/2}dx_n>0.
\end{align*}
We have a contradiction. Thus, $h_i$ is continuous on $C_\varepsilon$. In particular,
$$h_i\chi_{C_\varepsilon}\quad \text{ is $\mathcal{L}^{n-1}\res \Omega$-measurable for all $\varepsilon>0$.}$$
Since $\chi_{(\Omega\setminus C_\varepsilon)}\to 0$ in $L^1(\Omega)$, there exists a subsequence $\chi_{(\Omega\setminus C_{\varepsilon_k})}$ such that $\chi_{(\Omega\setminus C_{\varepsilon_k})}\to 0$ a.e. in $\Omega$. Therefore, $h_i\chi_{C_{\varepsilon_k}}\to h_i$ a.e. in $\Omega$ and hence $h_i$ is $\mathcal{L}^{n-1}\res \Omega\, $-measurable.

\noindent {\bf Step 3}: Now we define the following $\mathcal{L}^{n-1}\res\Omega\,$-measurable function
$$h(x')=\limsup_{i\to \infty}h_i(x').$$
Note that
\begin{align*}
\int_{\Omega}\left|m_E(f_i, x')-m_E(f, x')\right|dx'\leq 
\int_{\Omega}\int_{E_{x'}}|f_i(x',x_n)-f(x',x_n)|dx_ndx'\leq\|f_i-f\|_{L^1}\to 0.
\end{align*}
In particular, there exists a subsequence (still denoted as $f_i$) such that for $\mathcal{L}^{n-1}$-a.e. $x'$ in $\Omega$,
$$m_E(f_i, x')\to m_E(f, x').$$
Moreover,
\begin{align*}
&\left|\int_{-\infty}^{h_i(x')}f_i(x',x_n)dx_n-\int_{-\infty}^{h(x')}f(x',x_n)dx_n\right|\\
&=\left|\int_{-\infty}^{h_i(x')}f_i(x',x_n)dx_n-\int_{-\infty}^{h_i(x')}f(x',x_n)dx_n+\int_{-\infty}^{h_i(x')}f(x',x_n)dx_n-\int_{-\infty}^{h(x')}f(x',x_n)dx_n\right|\\
&=\left|\int_{-\infty}^{h_i(x')}f_i(x',x_n)-f(x',x_n)dx_n+\int_{h(x')}^{h_i(x')}f(x',x_n)dx_n\right|\\
&\leq \int_{\mathbb{R}}|f_i(x',x_n)-f(x',x_n)|dx_n+\left|\int_{h(x')}^{h_i(x')} f(x',x_n)dx_n\right|\to 0,
\end{align*}
where we have applied the dominated convergence theorem with $f(x',\cdot)\in L^1(\mathbb{R})$ and Step 1 (iii). Therefore,
\begin{align*}
m_E(f, x') =   \lim_{i \rightarrow \infty}  m_E(f_i, x') = \lim_{i \rightarrow \infty} \int_{-\infty}^{ h_i(x')}f_i(x', x_n)dx_n= \int_{-\infty}^{ h(x')}f(x', x_n)dx_n
\end{align*}
for $\mathcal{L}^{n-1}$-a.e. $x'$ in $\Omega$.
\qed

\subsection{Proof of Lemma \ref{l:meas implies meas arbitrary II}}
Let $\Omega$ be a bounded open set and $\delta>0$. We define the following $\mathcal{L}^{n-1}\res\Omega\,$-measurable function
$$m_{\delta}(f,x')=\min\{m_{E}(f,x')+\delta,m_{\mathbb{R}^n}(f,x')\}.$$
Let $f_i$ be the distribution function from Step 1 in the proof of Lemma \ref{l:meas implies meas arbitrary I}. Since
$$0\leq m_\delta(f_i,x')\leq \int_{-\infty}^{\infty}f_i(x',x_n)dx_n$$
and $f_i>0$, there exists a unique value $h^\delta_{i}(x')$ such that
$$m_\delta(f_i,x')=\int_{-\infty}^{h^\delta_i(x')}f_i(x',x_n)dx_n.$$
By using the same argument from the proof of Lemma \ref{l:meas implies meas arbitrary I}, we conclude that $h^\delta_i$ is $\mathcal{L}^{n-1}\res\Omega\,$-measurable. Define
$$h^\delta(x')=\limsup_{i\to \infty}h^\delta_i(x').$$
It is clear that $h^\delta$ is $\mathcal{L}^{n-1}\res\Omega\,$-measurable. Since $f_i\to f$ in $L^1(\mathbb{R}^{n})$,
$$m_E(f_i, x')\to m_E(f, x'),\quad m_{\mathbb{R}^n}(f_i, x')\to m_{\mathbb{R}^n}(f, x').$$
for $\mathcal{L}^{n-1}$-a.e. $x'$ in $\Omega$. In particular, we have $m_\delta(f_i, x')\to m_\delta(f, x')$ and
$$m_\delta (f,x')=\int_{-\infty}^{h^\delta(x')}f(x',x_n)dx_n$$
for $\mathcal{L}^{n-1}$-a.e. $x'$ in $\Omega$. Now we decompose $\Omega = Z \cup G$, where 
\begin{align*}
Z = \left\{x' : m_E(f,x') = m_{\mathbb{R}^n}(f,x')\right\}, \quad G= \Omega\setminus Z.
\end{align*}
Note that for $\mathcal{L}^{n-1}$-a.e. $x'$ in $\Omega$, the function $f(x', \cdot)$ is continuous and non-negative. Therefore, for $\mathcal{L}^{n-1}$-a.e. $x'$ in $G$, the collection $\{h^{\delta}(x')\}_{\delta}$ is monotonically decreasing as $\delta \rightarrow 0^+$. Let $\delta_k \rightarrow 0^+$ be a countable sequence, and define the function
\begin{align*}
    \tilde{h}(x') = \begin{cases}
        \lim\limits_{k \rightarrow \infty} h^{\delta_k}(x') & \mbox{$\mathcal{L}^{n-1}$-a.e. $x' \in G$}\\
       +\infty & x' \in Z. 
    \end{cases}
\end{align*}
We now claim that for $\mathcal{L}^{n-1}$-a.e. $x'$ in $\Omega$,
\begin{align}\label{largest_possible_func}
m_E(f,x')=\int_{-\infty}^{\tilde{h}(x')}f(x',x_n)dx_n.
\end{align}
By our definition of $\tilde{h}$, the equation \eqref{largest_possible_func} holds true for $x'\in Z$. For $\mathcal{L}^{n-1}$-a.e. $x' \in G$, we have $\tilde{h}(x')=\lim\limits_{k \rightarrow \infty} h^{\delta_k}(x')$ and 
$$m_E(f,x')\leq m_{\delta_k}(f,x')\leq m_E(f,x')+\delta_k.$$
Therefore,
\begin{align*}
    \left|m_E(f,x') - \int_{-\infty}^{ \tilde{h}(x')}f(x', x_n)dx_n\right|&\leq |m_E(f,x')-m_{\delta_k}(x')|+\left|\int_{\tilde{h}(x')}^{h^{\delta_k}(x')}f(x',x_n)dx_n\right|\to 0
\end{align*}
where we have applied the dominated convergence theorem with $f(x',\cdot)\in L^1(\mathbb{R})$. 

Next, we claim that $\tilde{h}$ is a largest possible value function of $m_E(f,\cdot)$.
Suppose there is another function $\alpha:\Omega\to [-\infty,\infty]$ such that 
\begin{align*}
    m_E(f,x') & = \int_{-\infty}^{ \alpha(x')}f(x', x_n)dx_n.
\end{align*}
Our goal is to show that $\alpha\leq \tilde{h}$ for $\mathcal{L}^{n-1}$-a.e. $x'$ in $\Omega$. It is easy to see that $\alpha(x')\leq \infty=\tilde{h}(x')$ for $x'\in Z$. For $\mathcal{L}^{n-1}$-a.e. $x' \in G$, we have $\tilde{h}(x')=\lim\limits_{k \rightarrow \infty} h^{\delta_k}(x')$, and there exists $N=N(x')\in \N$ such that for all $k\geq N$,
$$m_E(f,x')+\delta_k<m_{\mathbb{R}^n}(f,x').$$
Therefore, $m_{\delta_k}(f,x')=m_E(f,x')+\delta_k$ and
\begin{align*}
    \int_{-\infty}^{ \alpha(x')}f(x', x_n)dx_n=m_E(f,x')<m_{\delta_k}(f,x')=\int_{-\infty}^{h^{\delta_k}(x')}f(x',x_n)dx_n.
\end{align*}
The strict inequality above gives us
$$\alpha(x')<h^{\delta_k}(x')\quad \text{for all $k\geq N$}\implies \alpha(x')\leq \lim_{k\to \infty}h^{\delta_k}(x')=\tilde{h}(x').$$
Therefore, $\tilde{h}:\Omega\to [-\infty,\infty]$ is a largest possible value function of $m_E(f,\cdot)$. Finally, by the arbitrariness of $\Omega$, we may take an exhaustion of $\mathbb{R}^n$ by nested, bounded open sets and obtain by diagonalization a $\mathcal{L}^{n-1}$-measurable largest possible value function $h:\mathbb{R}^{n-1}\to [-\infty,\infty]$ of $m_E(f,\cdot)$.
\qed

\bibliography{references}
\bibliographystyle{amsalpha}

\end{document}